%% file: PFAVheuristics.tex
\documentclass{amsart}
\usepackage{amsmath}
\usepackage{amssymb}
\usepackage{rotating}
\usepackage{graphicx}
\usepackage{url}
\usepackage{array}
\usepackage{amscd}

\newtheorem{theorem}{Theorem}[section]
\newtheorem{refinedasymptoticformula}[theorem]{Estimate}
\newtheorem{fixedmaxrealsubfieldformula}[theorem]{Estimate}
\newtheorem{lemma}[theorem]{Lemma}
\newtheorem{proposition}[theorem]{Proposition}
\newtheorem{corollary}[theorem]{Corollary}

\theoremstyle{remark}
\newtheorem{remark}[theorem]{Remark}

\numberwithin{equation}{section}

\DeclareMathOperator\Gal{Gal}
\newcommand\case[1]{{\ensuremath{\mathrm(\mathit{#1}\mathrm)}}}

\newcommand\bZ{\ensuremath{\mathbb{Z}}}
\newcommand\bQ{\ensuremath{\mathbb{Q}}}
\newcommand\bR{\ensuremath{\mathbb{R}}}
\newcommand\bC{\ensuremath{\mathbb{C}}}
\newcommand\bF{\ensuremath{\mathbb{F}}}

\newcommand\cL{\ensuremath{\mathcal{L}}}
\newcommand\cM{\ensuremath{\mathcal{M}}}
\newcommand\cO{\ensuremath{\mathcal{O}}}

\newcommand\ve{\ensuremath{\varepsilon}}
\DeclareMathOperator\Tr{Tr}
\DeclareMathOperator\N{N}
\DeclareMathOperator\id{id}
\newcommand\gr{\ensuremath{\mathfrak{r}}}
\newcommand\ga{\ensuremath{\mathfrak{a}}}
\newcommand\gp{\ensuremath{\mathfrak{p}}}

\newcommand\gA{\ensuremath{\mathfrak{A}}}
\newcommand\gP{\ensuremath{\mathfrak{P}}}
\newcommand\gQ{\ensuremath{\mathfrak{Q}}}
\newcommand\abs[1]{\ensuremath{\lvert#1\rvert}}

\DeclareMathOperator\End{End}
\DeclareMathOperator\Aut{Aut}

\begin{document}
\title{Heuristics on pairing-friendly abelian varieties}

\date{\today}
\author[J. Boxall]{John Boxall}
\address{Laboratoire de Math\'ematiques Nicolas Oresme, CNRS -- UMR 6139, Universit\'e de Caen Basse-Normandie, boulevard mar\'echal Juin, BP 5186, 14032 Caen cedex, France}
\email{john.boxall@unicaen.fr}

\author[D. Gruenewald]{David Gruenewald}
\address{School of Mathematics and Statistics, University of Sydney, NSW 2006 Australia}
\email{davidg@maths.usyd.edu.au}

\begin{abstract}We discuss heuristic asymptotic formulae for the number of isogeny classes of pairing-friendly abelian varieties over prime fields, generalizing previous work of one of the authors \cite{Bo}.\end{abstract}



\thanks{This paper was written while the authors participated in the project \textit{Pairings and Advances in Cryptology for E-cash} (PACE) funded by the ANR}

\maketitle

\section*{Introduction}

Pairing-based cryptography uses non-degenerate pairings defined on a product $G_1\times G_2$ of two abelian groups and taking values in a third abelian group $G_3$. Typically, $G_1$, $G_2$ and $G_3$ are cyclic groups of the same prime order $r$. An important source of suitable groups is elliptic curves over finite fields, and in recent years generalizations using higher dimensional abelian varieties have been proposed.

As we shall recall briefly below, elliptic curves or abelian varieties possessing suitable subgroups for pairing-based cryptography satisfy very strong conditions, and are loosely referred to as pairing-friendly. These conditions suggest that they are very rare, and various upper bounds (either unconditional or depending on certain hypotheses) have been proposed, both for elliptic curves \cite{BK}, \cite{LS}, \cite{JULS} and more recently for Jacobians of genus two curves \cite{LaSh}. In \cite{Bo}, one of the authors investigated a heuristic asymptotic formula for the number of pairing-friendly elliptic curves over prime fields, and the purpose of the present paper is to present and provide computational evidence for generalizations of this to higher dimensional abelian varieties.

\subsection{Background} Let $q$ be a power of a prime $p$ and let $\bF_q$ denote a finite field with $q$ elements. Let $A$ be an abelian variety over $\bF_q$ of dimension $g\geq 1$. Let $\pi$ be the Frobenius endomorphism of $A$ over $\bF_q$, let $\ell\neq p$ be a prime and let $C_\pi$ be the characteristic polynomial of the action $\pi$ on the $\ell$-adic Tate module of $A$. Then Weil \cite{We} showed that $C_\pi$ is a monic polynomial of degree $2g$ with integer coefficients that is independent of $\ell$. Furthermore, he proved that the complex roots of $C_\pi$ are $q$-\emph{Weil numbers}, in other words algebraic integers of which any complex conjugate $\pi$ satisfies $\pi\bar{\pi}=q$. In general, by a $q$-\emph{Weil polynomial} we mean a monic polynomial with integer coefficients all of whose roots are $q$-Weil numbers. Thus, $C_\pi$ is a $q$-Weil polynomial of degree $2g$, and Tate \cite{Ta1} showed that $C_\pi$ depends only on the isogeny class of $A$. Furthermore, the order of the group $A(\bF_q)$ is equal to $C(1)$, so that the order of $A(\bF_q)$ is invariant under isogeny.  

A classification result was proved by Honda and Tate \cite{Ho}, \cite{Ta2}. They proved that there is a bijection between the set of $\bF_q$-isogeny classes of simple abelian varieties and the set of \emph{irreducible} Weil polynomials. Explicitly, the bijection associates to the isogeny class of the simple abelian variety $A$ the minimal polynomial $M_\pi$ of the Frobenius endomorphism $\pi$ of $A$. The characteristic polynomial $C_\pi$ is then a suitable power of $M_\pi$, so that $2\dim{A}$ is an even multiple of the degree of $M_\pi$. For a detailed study of this and related questions, we refer to \cite{Wa}. We simply mention the following result, which is part of Lemma~2.2 of \cite{FSS} and is proved using results in \cite{Wa}. 

\begin{proposition}\label{FSSprop} Let $g\geq 1$ be an integer, let $p$ be a prime and let $C$ be an irreducible $p$-Weil polynomial of degree $2g$. Then the simple abelian varieties over $\bF_p$ whose Frobenius endomorphism has minimal polynomial $C$ have dimension $g$. 
\end{proposition}

In what follows, if $g_0\geq 1$ is an integer, we understand by a \textit{triple of degree $g_0$} (or simply a triple if the reference to $g_0$ is clear) a triple $(r,C,q)$ where $r$ is a prime, $q$ a power of a prime and $C$ a $q$-Weil polynomial of degree $2g_0$ and $r$ divides $C(1)$. By our previous remarks, we can associate to (the isogeny class of) any pairing-friendly $g$-dimensional abelian variety $A$ over $\bF_q$ a triple $(r,C,q)$ of degree $g$, where $r$ is a prime such that $A(\bF_q)$ contains a subgroup $G_1$ of order $r$ and $C$ is the characteristic polynomial of the Frobenius endomorphism of $A$ over $\bF_q$. Conversely, if $q=p$ is prime and $C$ is irreducible of degree $2g$, Proposition~\ref{FSSprop} implies if $r$ is any prime dividing $C(1)$, then $(r,C,p)$ corresponds to an isogeny class of simple abelian varieties $A$ over $\bF_p$ such that $A(\bF_p)$ contains a subgroup $G_1$ of order $r$.

We now describe the conditions that a triple $(r,C,q)$ as above must satisfy assuming it is associated to an isogeny class of abelian varieties $A$ that are pairing-friendly. We refer to \cite{FST} (in the case of elliptic curves) and \cite{F} (in any dimension) for background and motivation. 

By definition, the \emph{rho-value} of the triple is $\rho=\rho(A)=\rho(r,C,q)=\frac{g\log{q}}{\log{r}}$, where $A$ is any member of the corresponding isogeny class and $g$ is the dimension of $A$.

(1)  Since $\#A(\bF_q)=C_\pi(1)$, it follows from the fact that the roots of $C_\pi$ are $q$-Weil numbers that $(\sqrt{q}-1)^{2g}\leq \#A(\bF_q)\leq (\sqrt{q}+1)^{2g}$. We deduce that, assuming $g$ is fixed and given any $\rho_0<1$, we must have $\rho>\rho_0$ whenever $q$ is sufficiently large.  In cryptographic applications, we want $\rho$ to be as close to $1$ as possible, since for fixed $r$, computations in the field $\bF_q$ will be faster. 

(2) There exists an integer $k\geq 2$ such that $r$ divides $\Phi_k(q)$. Here, $\Phi_k$ is the $k^{\text{th}}$ cyclotomic polynomial. (Some authors allow $k=1$, but we shall exclude this case. See Remark~\ref{kequalsone}.) Under some further mild restrictions, $A(\bF_{q^k})$ contains a subgroup $G_2$ of order $r$ different from $G_1$ and there is a computable non-degenerate pairing on $G_1\times G_2$ that takes values in the $r^{\text{th}}$ roots of unity in $\bF_{q^k}$ (see for example \cite{BK} when $g=1$ and Theorem 3.1 of \cite{RuSi}). We need to choose $k$ in such a way that the discrete logarithm problem in the field generated over the prime field $\bF_p$ by the $r^{\text{th}}$ roots of unity is  unfeasible, but not so large the computations in the field $\bF_{q^k}$ become unwieldy. We call $k$ the \emph{embedding degree} and $\frac{k}{g}$ the \emph{security ratio} of the triple (or of any abelian variety in the corresponding isogeny class). At present, the security ratio is in practice bounded by about $50$ when $g=1$ (see \cite{FST}). Although larger values of $k$ will be needed to maintain the same security level when $g$-dimensional abelian varieties are used, it is natural to consider $k$ as constant. 

Thus, when searching for abelian varieties of dimension $g$ for use in pairing-based cryptography, we seek triples $(r,C,q)$ with small embedding degree $k$ and rho-value as close to $1$ as possible. Once a suitable triple $(r,C,q)$ has been found, it is necessary to compute equations for some member of the corresponding isogeny class, the groups $G_1$ and $G_2$ and to be able to compute the pairing. The only known method of constructing such abelian varieties is the CM-method, which is based on the reduction at primes above $p$ of abelian varieties over number fields having complex multiplication, these being in turn constructed using theta functions or modular invariants. It is possible at present to construct elliptic curves with complex multiplication by the maximal order of an imaginary quadratic field whose discriminant is as large as about $10^{15}$ (see \cite{ES}).  On the other hand, only a few thousand explicit equations of genus two curves whose Jacobians have complex multiplication are known at present (see \cite{Ko}). Since this is computationally very heavy, it seems reasonable to consider triples that can be obtained from the reduction of a \emph{fixed} abelian variety, at least up to twist. 

\subsection{Presentation of the paper} Until now, when $g>1$, no examples of triples $(r,C,q)$ of cryptographic interest with rho-value less than or equal to $g$ have been found. The motivation that led to the present paper was an attempt to understand the reasons for this, and our main purpose is to discuss two heuristic estimates (see Estimates~\ref{Kfixedestimate} and \ref{Kplusfixedestimate}) related to the asymptotic growth as $x\to \infty$ for the number of triples $(r,C,p)$ of given degree $g$ with $r\leq x$, under certain conditions on the rho-value $g\frac{\log{p}}{\log{r}}$ and the polynomial $C$. 

In order to state these estimates formally, we need to consider triples of the form $(r,\pi,q)$, where $r$ is a prime and $\pi$ is a $q$-Weil number such that, if $K$ is any number field containing $\pi$, $r$ divides the $K/\bQ$-norm $\N_{K/\bQ}(\pi -1)$ of $\pi-1$.  This is because our heuristic arguments will be based mostly on the geometry of number fields.  See Remarks~\ref{autheuristics} and \ref{kplusauthheuristics} for how to interpret heuristics for triples $(r,\pi,q)$ in terms of triples $(r,C,q)$.
We say that $(r,\pi,q)$ has \textit{embedding degree} $k$ if $r$ divides $\Phi_k(q)$. If $2g$ is an even multiple of the degree of the number field $\bQ(\pi)$, we define the \textit{degree $g$ rho-value} of $(r,\pi,q)$ to be $g\frac{\log{q}}{\log{r}}$.  If $2g$ is equal to the degree of $\bQ(\pi)$, we simply call this the rho-value. 

From now on, we mostly restrict attention to triples $(r,\pi,p)$ with $p$ a prime. By Proposition~\ref{FSSprop}, such triples do correspond to isogeny classes of abelian varieties over prime fields. We write a triple as $(r,\pi, q)$ when we discuss matters in relation to arbitrary finite fields, it being understood that $q$ is a power of the prime $p$. 

To avoid repetition, we restrict attention to $g\neq 1$, referring the reader to \cite{Bo} and \cite{MS} for discussions of the elliptic curves case.

In both estimates we fix integers $g\geq 2$, $k\geq 2$ and a real number $\rho_0>1$ and suppose that the rho-value $g\frac{\log{p}}{\log{r}}$ is bounded by $\rho_0$. For a simple reason that will become clear in Lemma~\ref{trivialrhomin}, we in fact need to suppose that $\rho_0>\frac{g}{\varphi(k)}$, where $\varphi$ is Euler's totient function. In Estimate~\ref{Kfixedestimate}, we fix a CM-field $K$ of degree $2g$ and a CM-type $\Phi$ on $K$ and consider those triples $(r,\pi,p)$ of embedding degree $k$ with $\pi$ a $p$-Weil number belonging to $K$ that come from $\Phi$ (see \S~\ref{notandrev} for more explanation). In addition, we suppose that $\rho_0\neq g$. Then Estimate~\ref{Kfixedestimate} predicts that, except in a number of cases listed there, there is a constant $\alpha >0$ depending only on $k$ and $K$ such that the number $N(k,K,\Phi,\rho_0,x)$ of such triples $(r,\pi,p)$ with $r\leq x$ satisfies
\begin{equation}\label{introKfixed}
N(k,K,\Phi,\rho_0,x)\sim \frac{\alpha}{\rho_0}\int_{2}^{x}\frac{du}{u^{2-\frac{\rho_0}{g}}(\log{u})^2}, \qquad \text{as $x\to \infty$}.
\end{equation}
The constant $\alpha$ is explicit, and its value is implicit in the statement of Estimate~\ref{Kfixedestimate}. 
Since the integral on the right remains bounded as $x\to \infty$ when $\rho_0<g$, this asymptotic estimate would help to explain, if correct, the lack of known triples with small rho-values.

We have excluded the possibility $\rho_0=g$ since this is a borderline case, the integral being bounded as $x\to \infty$ when $\rho_0\leq g$ and unbounded when $\rho_0>g$.

%Although the integral remains bounded, in view of the way the rho-value is defined and the fact that $A(\bF_q)$ lies in the interval $[(\sqrt{q}-1)^{2g},(\sqrt{q}+1)^{2g}]$ it seems risky to conclude that $N(k,K,\Phi,g,x)$ is necessarily bounded when $x\to \infty$. Also, if $\rho_0>g$ we can integrate by parts to obtain 
%\begin{equation}\label{introKfixedbis}
%N(k,K,\Phi,\rho_0,x)\sim \frac{\alpha}{\rho_0(\rho_0-g)}\frac{x^{\frac{\rho_0}{g}-1}}{(\log{x})^2}, \qquad \text{as $x\to \infty$}.
%\end{equation}
%Since this is not justified when $\rho_0=g$ and we have no idea when $\rho_0>g$ which of the right hand sides of (\ref{introKfixed}) and (\ref{introKfixedbis}) should give a better numerical approximation to $N(k,K,\Phi,\rho_0,x)$, it seems risky to conclude that $N(k,K,\Phi,g,x)$ is necessarily bounded when $x\to \infty$.

Since our initial motivation was the search for abelian varieties with rho-value close to one, we were lead to study asymptotics for triples $(r,\pi,q)$ where the Weil numbers $\pi$ were allowed to  vary over larger sets. This lead to Estimate~\ref{Kplusfixedestimate}, where we fix a totally real field $K_0^+$ of degree $g$ and consider those $p$-Weil numbers that belong to some variable CM-field of degree $2g$ containing $K_0^+$. In addition, we suppose that $\rho_0\neq \frac{2g}{g+2}$. This estimate suggests that there is a constant $\beta>0$ depending only on $k$ and $K_0^+$ such that the number $R(k,K_0^+,\rho_0,x)$ of triples $(r,\pi,p)$ of embedding degree $k$ with $p$ prime, $r\leq x$ and $\pi$ a $p$-Weil number lying in some CM-field of degree $2g$ containing $K_0^+$ satisfies
\begin{equation*}
R(k,K_0^+,\rho_0,x)\sim \frac{\beta}{\rho_0}\int_2^x\frac{u^{\rho_0\big(\frac{1}{2}+\frac{1}{g}\big)-2}du}{(\log{u})^2} \qquad \text{as $x\to \infty$}.
\end{equation*} 
Again, the value of $\beta$ follows from the statement of Estimate~\ref{Kplusfixedestimate}. This estimate suggests that, if $g\leq 2$ and if $\rho_0>1$, then there should be infinitely triples $(r,\pi, p)$ with degree-$g$ rho value at most $\rho_0$. On the other hand, when $g\geq 3$, we only expect this to be thre case when $\rho_0>\frac{2g}{g+2}$.

We exclude the case $\rho_0=\frac{2g}{g+2}$ in Estimate~\ref{Kplusfixedestimate} for a similar reason that we exclude $\rho_0=g$ in Estimate~\ref{Kfixedestimate}.

We now outline briefly the contents of the paper. In \S~\ref{notandrev}, we recall some basic properties of CM-fields, CM-types and Weil numbers, and define the notion of a Weil number that comes from a given CM-type. In \S~\ref{Kfixedsection}, we state Estimate~\ref{Kfixedestimate} and give a heuristic argument in favour of it. In \S~\ref{Kfixedpolyfam}, we discuss the effect of polynomial families on the asymptotic formulae discussed in \S~\ref{Kfixedsection}. The discussion is analogous to that in \S~3 of \cite{Bo}. As in the case of elliptic curves, polynomial families having generic rho-value less than or equal to the dimension $g$ would produce counterexamples to the heuristics of \S~\ref{Kfixedpolyfam}. See Remark~\ref{counterexrem} for a discussion of the only known counterexample. Numerical computations in relation to Estimate~\ref{Kfixedestimate} are discussed in \S~\ref{Kfixedcomp}.  In \S~\ref{fixedmaxreal}, we state Estimate~\ref{Kplusfixedestimate} and present a heuristic argument that leads to it. The final \S~\ref{fixedmaxrealcomp} presents numerical computations in relation to the heuristics of \S~\ref{fixedmaxreal}.

\section{Notation and review of CM-types and Weil numbers}\label{notandrev}

\subsection{Notation and review of CM-types}\label{notandrevCM}

In this section we fix the embedding degree $k$ and the CM-field $K$ of degree $2g$ and we propose an asymptotic heuristic estimate as $x\to \infty$ for the number of triples $(r,\pi,p)$ using an approach similar to that already used for elliptic curves in \cite{Bo}. See also \S~3 of \cite{FSS}. Let $\zeta_k$ denote a primitive $k^{\text{th}}$ root of unity. If $F$ is a number field, we denote by $e(k,F)$ the degree of the number field $F\cap \bQ(\zeta_k)$. This is well-defined since $\bQ(\zeta_k)$ is a Galois extension of $\bQ$. Let $w_F$ denote the number of roots of unity in $F$ and $h_F$ the class number of $F$. We write $\N_{F/\bQ}$ for the absolute norm from $F$ to $\bQ$, applied to elements or to ideals. If $\alpha\in F$ and if $\sigma$ is an embedding of $F$ in another field or an automorphism of $F$, we write $\alpha^\sigma$ for the image of $\alpha$ under $\sigma$. If $F$ is a CM-field, we denote by $F^+$ the maximal real subfield of $F$ and by $c$ the non-trivial automorphism of $F$ that fixes $F^+$. We often write $\overline{\alpha}$ for the image of $\alpha\in F$ under $c$.

We briefly review the notion of CM-types, referring to \cite{Sh}, Chapter 2, for details. A \emph{CM-type} is a pair $(K,\Phi)$ (or simply $\Phi$ if the reference to $K$ is clear) consisting of a CM-field $K$ and a set $\Phi$ of $g$ embeddings $\phi:K\to \bC$ such that $\Phi\cup c\circ \Phi$ is the set of all embeddings of $K$ in $\bC$. Let $L$ denote a fixed Galois closure of $K$, so that $L$ is also a CM-field. Then, fixing an embedding of $L$ in $\bC$, we can view a CM-type $\Phi$ as a set of embeddings of $K$ in $L$. 

Let $(K,\Phi)$ and $(M, \Psi)$ two CM types. We say that $(K,\Phi)$ and $(M,\Psi)$ are \textit{equivalent} if there exists an isomorphism $\delta: K\to M$ and an automorphism $\iota$ of $\bC$ such that $\Psi=\{\delta^{-1}\circ  \phi\circ  \iota\mid \phi\in \Phi\}$.  If $M$ is also a subfield of the Galois closure $L$ of $K$ and we view CM-types on CM-subfields of $L$ as embeddings in $L$ as above, this definition is equivalent to asking that $\Psi=\{\delta^{-1}\circ  \phi\circ  \iota\mid \phi\in \Phi\}$ for some automorphisms $\delta$, $\iota$ of $L$. 

Write $G$ for the Galois group of $L$ over $\bQ$ and $H$ for the subgroup of $G$ corresponding to $K$. Any embedding $\phi$ of $K$ in $L$ can be extended to an automorphism of $L$, and if we also denote by $\phi$ one such extension, then all the extensions form a coset $\phi H$ of $H$. Let $S$ denote the set of all extensions of elements of $\Phi$ to automorphisms of $L$. Then $S$ is a CM-type on $L$ and the subgroup $H'=\{\gamma\in G\mid \gamma S=S\}$ contains $H$, and we say that $\Phi$ is \emph{primitive} if $H'=H$. 

Put $S^*=\{\sigma^{-1}\mid \sigma\in S\}$ and let $\hat{H}=\{\gamma\in G\mid \gamma S^*=S^*\}$. Then $\hat{H}$ is a subgroup of $G$ and the corresponding subfield $\hat{K}$ of $L$ is a CM-field, known as the \emph{reflex field} of $K$. Furthermore, the set $\hat{\Phi}$ of embeddings of $\hat{K}$ in $\bC$ obtained by restriction of elements of $S^*$ is a CM-type on $\hat{K}$, known as the \emph{reflex} of $\Phi$. Note that $\hat{\Phi}$ is always primitive.

By definition the \emph{$\Phi$-trace} is the map $\Tr_\Phi:K\to L$ that sends $\alpha\in K$ to $\sum_{\phi\in \Phi}\alpha^\phi$. One can show that the field $\hat{K}$ is generated over $\bQ$ by the set $\{\Tr_\Phi(\alpha)\mid \alpha\in K\}$, so that $\Tr_\Phi$ actually takes values in $\hat{K}$. It follows that the \emph{$\Phi$-norm} $\N_\Phi$, which is defined by $\N_\Phi(\alpha)=\prod_{\phi\in \Phi}\alpha^\phi$ for all $\alpha\in K$ also takes values in $\hat{K}$.  We define the \emph{$\hat{\Phi}$-trace} $\Tr_{\hat{\Phi}}$ and \emph{$\hat{\Phi}$-norm} $\N_{\hat{\Phi}}$ similarly; these maps are defined on $\hat{K}$ take values in the subfield $K'$ of $K$ corresponding to the subgroup $H'$ of $G$. The norm maps extend in an obvious way to maps on ideals $\N_\Phi: I_K\to I_{\hat{K}}$ and $\N_{\hat{\Phi}}: I_{\hat{K}}\to I_{K'}$, where $I_F$ denotes the group of fractional ideals of the number field $F$.

Let $Cl_{\hat{K}}$ be the ideal class group of $\hat{K}$ and and denote by $Cl(\hat{\Phi})$ the subgroup of $Cl_{\hat{K}}$ consisting of the ideal classes $\gamma$ such that $\N_{\hat{\Phi}}(\gamma)$ is the principal ideal class of $K$ and, if $\gA\in \gamma$, then the ideal $\N_{\hat{\Phi}}(\gA)$ of $K$ has a generator $\alpha$ such that $\alpha\overline{\alpha}$ is rational. (This makes sense as a class group since if $\gA$ is a principal ideal of $\hat{K}$ with generator $\beta$, then $\N_{\hat{\Phi}}(\beta)$ is a generator of $\N_{\hat{\Phi}}(\gA)$ and $\N_{\hat{\Phi}}(\beta)\overline{\N_{\hat{\Phi}}(\beta})=\N_{\hat{K}/\bQ}(\beta)$ is rational.) Let $h_{\hat{\Phi}}$ be the order of $Cl(\hat{\Phi})$.

\subsection{Weil numbers and characteristic polynomials}\label{Weilandchar}

Recall from the Introduction that if $q$ is a power of a prime $p$, a \textit{$q$-Weil number} is an algebraic integer all of whose complex conjugates $\pi$ satisfy $\pi\bar{\pi}=q$. 

Let $\pi$ be a $q$-Weil number, where $q=p^a$. If $\pi$ has a real conjugate, then $\pi^2=q$ and so $\pi=\pm{\sqrt{q}}$ and $\pi$ is totally real and belongs to $\bQ$ if $a$ is even and is real quadratic if $a$ is odd. On the other hand, if $\pi$ has no real conjugate, then $q/\pi$ is also a conjugate and $\pi$ is a root of the polynomial $X^2-\tau X+q$, where $\tau=\pi+q/\pi$ is a totally real algebraic integer. It follows that $\bQ(\pi)$ is a CM-field. Furthermore, every real conjugate of $\tau$ satisfies $\abs{\tau}\leq 2\sqrt{q}$. Conversely, if $\tau$ is a totally real algebraic integer all of whose real conjugates satisfy this inequality, then the roots of $X^2-\tau X+q$ are $q$-Weil numbers. 

Recall that an abelian variety $A$ over $\bF_q$ of dimension $g$ is said to be \textit{ordinary} if the group $A[p](\overline{\bF}_q)$ of $p$-torsion points over an algebraic closure $\overline{\bF}_q$ of $\bF_q$ has $p$-rank $g$. We refer to \cite{Wa}, \S~7 for the following result.

\begin{proposition}\label{TFAEordinary} Let $A$ be a simple abelian variety over $\bF_q$ of dimension $g$ and let $\pi$ be the Frobenius endomorphism of $A$ and let $\tau=\pi+q/\pi$. Then the following are equivalent.

\case{i} $A$ is ordinary.

\case{ii} $\pi$ and $q/\pi$ are coprime algebraic integers.

\case{iii} $\tau$ and $q$ are coprime algebraic integers.

Furthermore, if these conditions are satisfied, $\End(A)\otimes \bQ=\bQ(\pi)$ is a CM-field. 
\end{proposition} 

The fact that $\bQ(\pi)$ is a CM-field follows from the previous discussion. We deduce the following Corollary to Proposition~\ref{FSSprop}, already remarked in \cite{FSS}.

\begin{corollary} Let $g\geq 1$ be an integer and let $K$ be a CM-field of degree $2g$. Let $p$ be a prime and let $C$ be an irreducible $p$-Weil polynomial of degree $2g$ such that $K\simeq \bQ[X]/C(X)$. If $p$ is unramified in $K$, then the abelian varieties over $\bF_p$ belonging to the isogeny class corresponding to $C$ are ordinary. 
\end{corollary}

\begin{proof}
Indeed, if $p$ is such that the abelian varieties are not ordinary, then $\pi$ and $p/\pi$ have a common prime ideal factor $\gp$. But then $\gp^2$ divides $p$, so $p$ is ramified in $K$.
\end{proof}

Inspired by this result, we say that the $q$-Weil number $\pi$ is \textit{ordinary} if it satisfies the equivalent conditions \case{ii} and \case{iii} of the Proposition. 

Let $\pi$ be a $q$-Weil number belonging to the CM-field $K$ of degree $2g$. Then by the \textit{characteristic polynomial of $\pi$}, denoted by $C_\pi$, we mean the characterstic polynomial of the endomorphism multiplication-by-$\pi$ of the $\bQ$-vector space $K$. Then $C_\pi$ is a power of the minimal polynomial of $\pi$ and depends only on $g$ and not on $K$. In particular, if $\pi'$ is another $q$-Weil number belonging to $K$, then $C_{\pi'}=C_\pi$ if and only if $\pi'$ is a conjugate of $\pi$. Furthermore, if $K=\bQ(\pi)$, then $C_\pi$ is the minimal polynomial of $\pi$ and, if  $\Aut(K)$ denotes the automorphism group of $K$, then the number of conjugates of an element of $K$ is equal to the order of $\Aut(K)$. Thus $\# \Aut(K)$ triples of the form $(r,\pi,q)$ give rise to the same triple $(r,C,q)$.

\subsection{Weil numbers coming from a given CM-type}
 
Let $(K,\Phi)$ be a CM-type. If $\pi\in K$ is a Weil number, we say that $\pi$ \emph{comes from $\Phi$} if there is an ideal $\gA\in I_{\hat{K}}$ such that $\N_{\hat{\Phi}}(\gA)$ is a principal ideal of $K$, generated by $\pi$. Similarly, we say that the triple of one of the forms $(r,C,q)$ and $(r,\pi, q)$ as above comes from $\Phi$ if $\pi$ comes from $\Phi$.  

Before turning to heuristic arguments, we derive some unconditional results about Weil numbers coming from a given CM-type.   

\begin{proposition}\label{primitiveweil}
Let $(K,\Phi)$ be a CM-type, let $p$ be a prime unramified in $K$ and let $\pi\in K$ be a $p$-Weil number coming from $\Phi$.

\case{i} There is a unique prime ideal $\gP$ of $\hat{K}$ such that $\pi$ generates the ideal $\N_{\hat{\Phi}}\gP$ of $K$. Furthermore, $\gP$ is of degree one, and its ideal class belongs to $Cl(\hat{\Phi})$.

\case{ii} If $(K,\Phi)$ is primitive, then $K=\bQ(\pi)$.
\end{proposition}

\begin{proof} \case{i} By considering the factorization of $p$ as a product of prime ideals in $L$ and $K$, and using $\pi\overline{\pi}=p$, one sees that there is an integral ideal $\gP$ of $\hat{K}$ such that $\N_{\hat{\Phi}}\gP$ is equal to the ideal generated by $\pi$. Furthermore, one necessarily has $\N_{\hat{K}/\bQ}(\gP)=p$, so that since $p$ is prime $\gP$ is necessarily a prime ideal of degree one. It is clear from the definitions that the ideal class of $\gP$ must belong to $Cl(\hat{\Phi})$.

Suppose that $\gP'$ is a second prime ideal of $\hat{K}$ such that $\N_{\hat{\Phi}}\gP'$ is equal to the ideal generated by $\pi$. Let $\gQ$ and $\gQ'$ denote respectively prime ideals of $L$ dividing $\gP$ and $\gP'$. Since $p$ divides both $\gQ$ and $\gQ'$, there is an element $\tau \in G$ such that $\gQ'=\gQ^\tau$. Let $G_\gQ$ be the decomposition group of $\gQ$; since $\gP$ has degree one, $G_\gQ$ is a subgroup of $\hat{H}$. Since $L$ is the Galois closure of $K$, $p$ is unramified in $L$ and the ideal factorization of $\gP$ in $L$ is $\prod_{\sigma\in G_\gQ\backslash\hat{H}}\gQ^\sigma$. From the definition of $\N_{\hat{\Phi}}$, we see that the ideal factorization of $\pi$ in $L$ is equal to $\prod_{\sigma\in G_\gQ\backslash S^*}\gQ^\sigma$. By hypothesis, it is also equal to $\prod_{\sigma\in G_{\gP^\tau}\backslash S^*}\gQ^{\tau\sigma}$. It follows that $\tau S^*=S^*$, so that $\tau\in \hat{H}$. But then $\gQ$ and $\gQ'$ divide the same prime ideal of $\hat{K}$, so that $\gP'=\gP$.

\case{ii} It suffices to prove that if $\gamma\in G$ is such that $\pi^\gamma$ and $\pi$ have the same ideal factorization in $L$, then $\gamma\in H$. To do this, we return to the ideal factorization 
$\prod_{\sigma\in G_\gQ\backslash S^*}\gQ^\sigma$ of $\pi$. The ideal factorization of $\pi^\gamma$ is then $\prod_{\sigma\in G_\gQ\backslash S^*}\gQ^{\sigma\gamma}$, and this can only coincide with that of $\pi$ if $S^*\gamma=S^*$, or equivalently if $\gamma S=S$. Since $\Phi$ is primitive, this implies that $\gamma\in H$.\end{proof}

Despite its simple proof, we have been unable to find a reference to the following result in the literature.

\begin{theorem}\label{asymptoticWeil}
Let $\Phi$ be a CM-type on $K$. Then the number $\pi_\Phi(x)$ of $p$-Weil numbers coming from $\Phi$ with $p$ prime and $p\leq x$ is asymptotically equal to
\begin{equation}
\pi_\Phi(x)\sim \frac{w_Kh_{\hat{\Phi}}}{h_{\hat{K}}}\int_2^x\frac{du}{\log{u}}
\end{equation}
as $x\to \infty$. 
\end{theorem}

\begin{proof} Let $\pi$ be a $p$-Weil number coming from $\Phi$. Since only finitely many primes ramify in $K$, we can suppose $p$ unramified. By Proposition~\ref{primitiveweil}~\case{i}, there is a unique prime ideal of degree one $\gP$ of $\hat{K}$ such that $\N_{\hat{\Phi}}(\gP)$ is equal to the ideal of $K$ generated by $\pi$. On the other hand, by Kronecker's lemma on roots of unity, every other Weil number that generates the same ideal as $\pi$ is of the form $\zeta \pi$ for some root of unity $\zeta\in K$. It follows that $\pi_\Phi(x)$ is equal to $w_K$ times the number of degree one prime ideals $\gP$ of $\hat{K}$ with $\N_{\hat{K}/\bQ}(\gP)\leq x$. Since the ideal class of $\gP$ must belong to $Cl(\hat{\Phi})$, the result thus follows by applying the prime ideal theorem in number fields (see for example \cite{Na}~Chapter VII).\end{proof}

\section{A heuristic asymptotic formula for a fixed CM-type} \label{Kfixedsection}
Before discussing Estimate~\ref{Kfixedestimate}, we begin with the following simple observation. As before, $\varphi$ denotes Euler's totient function. Also, we allow triples $(r,\pi,q)$ with $q$ a prime power.

\begin{lemma}\label{trivialrhomin}
Let $g$ and $k\geq 1$ be fixed. If $\rho_0< \frac{g}{\varphi(k)}$, then there are only finitely many triples $(r,\pi,q)$ of genus $g$ with embedding degree $k$ and $g$-degree rho-value $\leq \rho_0$. 
\end{lemma}

\begin{proof}  Suppose that there are infinitely many triples $(r,\pi,q)$ with embedding degree $k$ and rho-value $\rho\leq \rho_0$. For any fixed $r$, since $q\leq r^{\frac{\rho_0}{g}}$, there are only finitely many possibilities for the prime-power $q$ such that $(r,\pi,q)$ is a triple whose first member is $r$. On the other hand, by considering the ideal factorization of $q$ in $K$ and using Kronecker's lemma on roots of unity, we deduce that for fixed $q$ there are only finitely many $q$-Weil numbers belonging to $K$. Hence the set of values $r$ that appear as first members of triples is unbounded. Recall that the $k^{\text{th}}$ cyclotomic polynomial $\Phi_k$ is a monic polynomial of degree $\varphi(k)$. Since $r$ divides $\Phi_k(q)$, we have $r\leq \Phi_k(q)$. On the other hand, $\Phi_k(q)\sim q^{\varphi(k)}$ and $q^{\varphi(k)}\leq r^{\frac{\varphi(k)\rho}{g}}$, so that $\Phi_k(q)\leq 2r^{\frac{\varphi(k)\rho}{g}}\leq 2r^{\frac{\varphi(k)\rho_0}{g}}$ if $r$ is large enough. We deduce that $r\leq 2r^{\frac{\varphi(k)\rho_0}{g}}$ if $r$ is large enough. Since $\rho_0<\frac{g}{\varphi(k)}$, this is impossible.\end{proof}

The purpose of rest of this section is to give a heuristic argument in support of the following asymptotic estimate.

\begin{refinedasymptoticformula}  \label{Kfixedestimate}
Let $g\geq 2$, $k\geq 2$ be integers, and let $\rho_0>\max(1,\frac{g}{\varphi(k)})$ be a real number such that $\rho_0\neq g$. Fix a CM-field $K$ of degree $2g$, a CM-type $\Phi$ on $K$ and let $e(k,K)$, $w_K$, $h_{\hat{K}}$ and $h_{\hat{\Phi}}$ be as above.  Then the number of triples $(r,\pi,p)$ as above with $r\leq x$ and $p\leq r^{\frac{\rho_0}{g}}$ that come from $\Phi$ is equivalent as $x\to \infty$ to   
\begin{equation*}
\frac{e(k,K)gw_Kh_{\hat{\Phi}}}{\rho_0h_{\hat{K}}}\int_{2}^x\frac{du}{u^{2-\frac{\rho_0}{g}}(\log{u})^2}
\end{equation*}
\end{refinedasymptoticformula}

\begin{remark} \label{autheuristics}  It follows from Proposition~\ref{primitiveweil}~\case{ii}  and the discussion at the end of \ref{Weilandchar} that, when $\Phi$ is primitive, the number of triples $(r,C,p)$ with $K\simeq \bQ[x]/C(x)$ is expected to be asymptotically equivalent to
\begin{equation*}
\frac{e(k,K)gw_Kh_{\hat{\Phi}}}{\#\big(\Aut(K)\big)\rho_0h_{\hat{K}}}\int_{2}^x\frac{du}{u^{2-\frac{\rho_0}{g}}(\log{u})^2}
\end{equation*}
as $x\to \infty$. Also, when $p$ is unramified in $K$, the triple is ordinary. 
\end{remark}

We now complete the heuristic argument which will lead to Estimate~\ref{Kfixedestimate}, in a manner similar to \cite{Bo}.

Let $r$ be given. The probability that $r$ is a prime and splits completely in $\bQ(\zeta_k)$ is equal to the probability that $r$ is prime and that $r\equiv 1\pmod{k}$ which is $\frac{1}{\varphi(k)\log{r}}$. On the other hand, when $r\equiv 1\pmod{k}$, $\Phi_k$ has $\varphi(k)$ distinct roots $\pmod{r}$. Thus if $p$ is any integer, the probability that $\Phi_k(p)\equiv 0\pmod{r}$ is roughly $\frac{\varphi(k)}{r}$. Hence the probability $r$ is prime and divides $\Phi_k(p)$ is roughly $\frac{1}{\varphi(k)\log{r}}\frac{\varphi(k)}{r}=\frac{1}{r\log{r}}$.

On the other hand, we want $p$ to be a prime and $\pi\in K$ to be a $p$-Weil number such that $r$ divides $\N_{K/\bQ}(\pi-1)$. We ignore the case where $r^2$ divides $\N_{K/\bQ}(\pi-1)$, assuming that as $r$ increases it occurs with negligible frequency. This means that there is a unique degree one prime ideal $\gr$ dividing $r$ that also divides $\pi-1$.  We assume that the Weil numbers $\pi$ behave randomly with respect to division by a non-zero ideal of $K$. This means that the probability that $\pi-1$ is divisible by $\gr$ is $\frac{1}{r}$.  In any number field, the average number of primes of degree one dividing a given rational prime is one, so that this is equivalent to the probability that $r$ divides $\N_{K/\bQ}(\pi-1)$ being equal to $\frac{1}{r}$. 

Now  the events that $r$ splits completely in $\bQ(\zeta_k)$ and there exists a degree one prime in $K$ dividing $r$ are not independent in general. This is the case if and only if $K$ and $\bQ(\zeta_k)$ are linearly disjoint over $\bQ$. In general, the probability that $r$ is a prime that enjoys both of these properties is $\frac{e(k,K)}{r^2\log{r}}$.

Finally, from Theorem~\ref{asymptoticWeil}, we see that the expected number of $p$-Weil numbers coming from $\Phi$ with prime $p \leq r^{\frac{\rho_0}{g}}$ is about $\frac{w_Kh_{\hat{\Phi}}}{h_{\hat{K}}}\frac{r^{\frac{\rho_0}{g}}}{\log{r^{\frac{\rho_0}{g}}}}$, so that the total number of triples $(r,\pi,p)$ satisfying the hypotheses of Estimate~\ref{Kfixedestimate} is expected to be asymptotically equivalent to
\begin{equation*}
\sum_{2\leq r\leq x}\frac{e(k,K)}{r^2(\log{r})}\frac{w_Kh_{\hat{\Phi}}}{h_{\hat{K}}}\frac{r^{\frac{\rho_0}{g}}}{\log{r^{\frac{\rho_0}{g}}}}=\frac{e(k,K)gw_Kh_{\hat{\Phi}}}{\rho_0h_{\hat{K}}}\sum_{2\leq r\leq x}\frac{1}{r^{2-\frac{\rho_0}{g}}(\log{r})^2},
\end{equation*}
where the sums are over all integers $r$ between $2$ and $x$. Replacing the sum by an integral leads to Estimate~\ref{Kfixedestimate}. 

\begin{remark} \label{kequalsone} 
Suppose that $\bF_q$ contains the $r^{\text{th}}$ roots of unity. Let $A/\bF_q$ be a simple ordinary abelian variety of dimension $g$ and suppose that the prime $r\neq p$ divides the order of $A(\bF_q)$. Then $q\equiv 1\pmod{r}$ and hence $1\pmod{r}$ is a root of multiplicity at least $2$ of the characteristic polynomial of the Frobenius endomorphism $\pi$ of $A$ acting on the $\bF_r$-vector space $A[r]$ of points of order dividing $r$ of $A$. Suppose that $r$ is prime to the discriminant of the order $\bZ[\pi,p/\pi]$, which is contained in $\End(A)$. Then $r$ can be factored as a product of distinct proper $\End(A)$-ideals, say $r\End(A)=\gr_1\gr_2\dots \gr_h$. We can then write $A[r]$ as direct sum
\begin{equation*}
A[r]=A[\gr_1]\oplus A[\gr_2]\oplus \cdots \oplus A[\gr_h],
\end{equation*}
where for each $i$, $A[\gr_i]$ denotes the subgroup of $A[r]$ killed by all the elements of $\gr_i$. We deduce that the action of $\pi$ on $A[r]$ is semi-simple, so that the $1$-eigenspace of $\pi$ in $A[r]$ viewed as a $\bF_r$-vector space is of dimension at least two. It follows that $r^2$ necessarily divides the order of $A(\bF_q)$. Thus the assumption made in the study of Estimate~\ref{Kfixedestimate} that we can ignore cases where $r^2$ divides the order $A(\bF_q)$ is not justified when $k=1$.  
\end{remark}

\section{The influence of complete polynomial families}\label{Kfixedpolyfam}

As in the genus one case (see \S~3 of \cite{Bo}), complete polynomial families with small generic rho-value can be expected to produce more triples than predicted by Estimate~\ref{Kfixedestimate}. The argument is a simple generalization of the proof of Theorem~3.1 of \cite{Bo}. For a detailed discussion of the construction of such families in genus $2$ or more, we refer to \cite{F}. Examples of such families occur in \cite{F}, \cite{FSa}, \cite{Ka}, \cite{KT}, \cite{D} and \cite{GuVe}. See also \cite{GMV}. We recall briefly the definitions given in \cite{F}.

As in the previous section, we fix $k$ and the CM-field of $K$ of degree $2g$. We fix a CM-type $\Phi$ on $K$ and denote by $\hat{K}$ the reflex field. Let $w$ be an indeterminate. We denote by $\overline{\alpha}(w)$ the polynomial in $K[w]$ obtained from $\alpha(w)$ by applying complex conjugation to the coefficients. If $\alpha(w)\in K[w]$, then $\N_{K/\bQ}(\alpha(w))=\prod_{i=1}^{2g}(\alpha(w))^{\sigma_i}$ where, for each $i$, $(\alpha(w))^{\sigma_i}$ is the polynomial obtained by applying the complex embedding $\sigma_i$ of $K$ to each coefficient of $\alpha(w)$. We extend the definitions of $\N_{\Phi}$ and $\N_{\hat{\Phi}}$ to polynomials in a similar manner. Following \cite{F} and generalizing \cite{BWe} and \cite{FST}, a \emph{polynomial family} (or complete family) is a pair of polynomials $(r_0(w),\pi_0(w))$ with $r_0(w)\in \bQ[w]$ and $\pi_0(w)\in K[w]$ such that:

(\textit{i}) $r_0(w)$ is irreducible and the field $\bQ[w]/r_0(w)$ contains a subfield isomorphic to $\hat{K}$.

(\textit{ii}) $\pi_0(w)\overline{\pi}_0(w)$ is a polynomial $p_0(w)\in \bQ[w]$,

(\textit{iii}) $r_0(w)$ divides both $\Phi_k(p_0(w))$ and $\N_{K/\bQ}(\pi_0(w)-1)$, 

(\textit{iv}) there exist infinitely many integers $w_0$ such that $r_0(w_0)$ is a prime (or a near-prime) and $p_0(w_0)$ is a prime.  

When $w_0$ is chosen such that $r_0(w_0)$ and $p_0(w_0)$ are prime, then $\pi_0(w_0)$ is a $p_0(w_0)$-Weil number and the characteristic polynomial of $\pi_0(w_0)$ is $C_{\pi_0(w_0)}(x)=\N_{K/\bQ}(x-\pi_0(w_0))$. Furthermore, the condition (\textit{iii}) implies that $r_0(w_0)$ divides both $\Phi_k(p_0(w_0))$ and $C_{\pi_0(w_0)}(1)$, and therefore give rise to a triple $(r_0(w_0),\pi_0(w_0),p_0(w_0))$.

As $w_0$ tends to $\infty$, the rho-value of the triple tends to $\frac{g\deg p_0}{\deg r_0}$, so we call this quantity the \emph{generic rho-value} of the family. 

Following work of Hardy and Littlewood, Schinzel and others, Bateman and Horn \cite{BH} (see also Conrad \cite{K}) give heuristic asymptotic formulae concerning sets of polynomials taking prime values simultaneously. Let $f_1$, $f_2$, \dots, $f_m\in \bQ[w]$ be irreducible and suppose that $f_i(\bZ)\subseteq \bZ$ for all $i\in \{1,2,\dots, m\}$. These heuristics imply that, as $X\to \infty$, the number of integers $w_0$ with $0\leq w_0\leq X$ such that $f_i(w_0)$ is a prime number for all $i$ is either finite or asymptotically equivalent to
\begin{equation} \label{BHasymptotic}
A\int_{2}^X\frac{du}{(\log{u})^m}\sim A\frac{X}{(\log{X})^m} \quad \text{as} \quad X\to \infty,
\end{equation}
where the constant $A> 0$ depends only on the polynomials $f_i$ and is given explicitly in \cite{BH} (see also the interesting discussion in \cite{K}). 

\begin{theorem} Keep the notations used above and assume the Bateman-Horn heuristics (\cite{BH}, \cite{K}). 

\case{i} Suppose that the polynomial family parametrized by $r_0(w)$ and $\pi_0(w)$ produces, asymptotically as $x\to \infty$, more triples $(r,\pi,p)$ with $r\leq x$ than predicted by Estimate~\ref{Kfixedestimate}. Then $g\frac{\deg{p_0}}{\deg{r_0}}\leq \rho_0\leq g(1+\frac{1}{\deg{r_0}})$. Furthermore, the generic rho-value of the family is less than or equal to $g$.

\case{ii} Conversely, suppose that $g\frac{\deg{p_0}}{\deg{r_0}}<\rho_0< g(1+\frac{1}{\deg{r_0}})$.  Then the polynomial family produces more triples than predicted by Estimate~\ref{Kfixedestimate}.
\end{theorem} 

\begin{proof}
The proof is the same as that of Theorem~3.1 of \cite{Bo}, so we only sketch the details.

\case{i} Using (\ref{BHasymptotic}) with $m=2$, we find the number of triples $(r,\pi,p)$ with $r\leq x$ that the family is expected to generate is asymptotically equivalent to 
\begin{equation} \label{BTKheuristic}
C\frac{x^{1/\deg{r_0}}}{(\log{x})^2} \quad \text{ as $x\to \infty$},
\end{equation}
where $C>0$ is another constant.

On the other hand, Estimate~\ref{Kfixedestimate} predicts that the number of triples $(r,\pi,p)$ with $r\leq x$ should be finite when $\rho_0< g$ and equivalent to
\begin{equation} \label{Kfixedestimateip}
\frac{e(k,K)g^2w_Kh_{\hat{\Phi}}}{\rho_0(\rho_0-g)h_{\hat{K}}}\frac{x^{\frac{\rho_0}{g}-1}}{(\log{x})^2}
\end{equation}
as $x\to \infty$ when $\rho_0>g$. In any case, we see that the complete family is expected to provide more examples than predicted by Estimate~\ref{Kfixedestimate} when $\frac{1}{\deg{r_0}}>\frac{\rho_0}{g}-1$, but not when $\frac{1}{\deg{r_0}}<\frac{\rho_0}{g}-1$. 

Now, as $|w_0|\to \infty$, the rho-value of the triple $(r_0(w_0),\pi_0(w_0),p_0(w_0))$ converges to the generic rho-value of the polynomial family, which is $g\frac{\deg{p_0}}{\deg{r_0}}$. Thus, since by hypothesis the family produces more examples than predicted by Estimate~\ref{Kfixedestimate}, we have
\begin{equation*}
g\frac{\deg{p_0}}{\deg{r_0}}\leq \rho_0\leq g\big(1+\frac{1}{\deg{r_0}}\big).
\end{equation*}

On the other hand, points \case{i} and \case{ii} of the definition of a complete family show that $\deg{r_0}$ and $\deg{p_0}$ are even integers, so we deduce that $\deg{p_0}\leq \deg{r_0}$ which in turn implies that the generic rho-value of the family is at most $g$. This proves \case{i}.

\case{ii} Conversely, if $\frac{g\deg{p_0}}{\deg{r_0}}<\rho_0<g(1+\frac{1}{\deg{r_0}})$ then, when $\abs{w_0}$ is sufficiently large, the rho-value of $(r_0(w_0),\pi_0(w_0),p_0(w_0))$ is less than $\rho_0$, and reversing the argument used to prove the first assertion of \case{i} shows that the family contains more triples than predicted by Estimate~\ref{Kfixedestimate}.   
\end{proof}

\begin{remark} 
What happens when $\rho_0=g(1+\frac{1}{\deg{r_0}})$ depends on the relative sizes of the constants $C$ and $\frac{e(k,K)g^2w_Kh_{\hat{\Phi}}}{\rho_0(\rho_0-g)h_{\hat{K}}}$  appearing in (\ref{BTKheuristic}) and in (\ref{Kfixedestimateip}). We do not know what happens when $\rho_0=g\frac{\deg{p_0}}{\deg{r_0}}$.
\end{remark}

\begin{remark} \label{counterexrem}
The inequality $\deg{p_0}\leq \deg{r_0}$ in \case{i} of the theorem implies that if the family produces more examples than predicted by Estimate~\ref{Kfixedestimate}, then its generic rho-value is at most $g$. When $g=1$, the only known polynomial family that satisfies this condition and that produces curves over prime fields is the well-known Barreto-Naehrig family of prime-order elliptic curves, where 
\begin{equation*} 
r_0(w)=36w^4+36w^3+18w^2+6w+1,\qquad  \pi_0(w)=\frac{t_0(w)+\sqrt{-3}y_0(w)}{2},
\end{equation*}
with $t_0(w)=6w^2+1$ and $y_0(w)=6w^2+4w+1$. Here $k=12$ and $K=\bQ(\sqrt{-3})$, and the family provides a genuine counterexample to Estimate~\ref{Kfixedestimate} when $\rho_0<1.25$. This is discussed in detail in \cite{Bo}. 

By taking abelian varieties isogenous to the $g^{\text{th}}$ power of the elliptic curves produced by the Barreto-Naehrig family, we obtain for every $g\geq 1$ a polynomial family of abelian varieties with $k=12$, $K$ any CM-field containing $\bQ(\sqrt{-3})$, and generic rho-value $g$. For the CM-type on $K$ we take any extension of a CM-type on $\bQ(\sqrt{-3})$. These families produce more triples than predicted by Estimate~\ref{Kfixedestimate} when $\rho_0<\frac{5g}{4}$.     
\end{remark}

\begin{remark} \label{genrhovalues}
In \cite{F}, Freeman gives a general construction of Brezing-Weng families that is inspired by the construction of Weil numbers as images under the $\hat{K}$-norm of elements of $\hat{K}$. This gives families whose generic rho-value is usually a little smaller than $2g\hat{g}$, where $2\hat{g}=[\hat{K}:\bQ]$. However he does give an example with $k=5$, $K=\bQ(\zeta_5)$ and generic rho-value $4$. When $g=2$, a somewhat different approach is discussed in \cite{D} (see also \cite{LaSh}), but when applied to primitive CM-types again gives rho-values close to $8$. To the best of our knowledge, no example of a family of absolutely simple abelian surfaces (which correspond to primitive CM-types) with generic rho-value less than $4$ is known. Apart from the case $g=2$ of the generalized Barreto-Naehrig construction mentioned in Remark~\ref{counterexrem}, many examples of polynomial families of abelian surfaces that are not absolutely simple with generic rho-value less than 4 are known (\cite{FSa}, \cite{Ka}, \cite{KT}, \cite{D}, \cite{GuVe}), but in every case have a generic rho-value strictly greater than $2$. 

Thus, while it seems difficult to construct polynomial families with generic rho-value $g$ or less, we see no reason why the Barreto-Naehrig family should be the only such family.
\end{remark}

\section{Numerical evidence in the fixed CM-field case}  \label{Kfixedcomp}

 In this section, we report on numerical evidence for Estimate~\ref{Kfixedestimate}.

\subsection{Examples in genus $2$ and $3$}

% In this section, we report on numerical evidence for Estimate~\ref{Kfixedestimate} when $g=2$ and $g=3$.  

For the convenience of the reader we summarize briefly all possible CM-types up to equivalence for $g=2$ and $3$. Recall that $K$ denotes a CM-field of degree $2g$, $K^+$ the maximal real subfield of $K$, $L$ a fixed Galois closure of $K$ and $G$ the Galois group of $L$ over $\bQ$. If $F$ is a subfield of $L$, we denote by $G^F$ the subgroup of $G$ that fixes $F$. We fix an embedding of $L$ into $\bC$. When we write a CM-type $\Phi$ on $K$ as a set $\{\tau_1,\tau_2,\dots \}$ of elements of $G$, we mean that $\Phi$ is the restriction to $K$ of these elements.  Similarly for the reflex type $\hat{\Phi}$.

\noindent\underline{$g=2$} (cf. \cite{Sh}, pages 64--65). There are three possibilities for $K$.

\case{i} $K=L$ is a Galois biquadratic extension of $\bQ$. Let $K_0$ be an imaginary quadratic subfield of $K$ and let $\tau$ be the non-trivial element of $G^{K_0}$. Then $\Phi=\{\id, \tau\}$ is a CM-type on $K$ extending the CM-type $\{\id\}$ of $K_0$.  We have $\hat{K}=K_0$ and $\hat{\Phi}=\{\id\}$. Hence there are two equivalence classes of CM-types on $K$ corresponding to the two imaginary quadratic subfields. Neither of them is primitive.

\case{ii} $K=L$ is a Galois cyclic extension of $\bQ$. If $\tau$ is a generator of $G$, then every CM-type is equivalent to $\Phi=\{\id, \tau\}$. This is a primitive CM-type, $\hat{K}=K$ and $\hat{\Phi}=\{\id,\tau^{-1}\}$.

\case{iii} $K$ is not Galois over $\bQ$, and $G$ is a dihedral group generated by $\sigma$ and $\tau$ with $\sigma$ of order $2$ having $K$ as fixed field and $\tau$ is of order $4$. Every CM-type is equivalent to  $\Phi=\{\id, \tau\}$. Again this is a primitive CM-type, $\hat{K}$ is the field fixed by $\tau\sigma$ and $\hat{\Phi}=\{\id, \sigma\}$. 

\noindent\underline{$g=3$}. There are four possibilities for $K$.

\case{iv} $K=L$ is a degree $6$ Galois cyclic extension of $\bQ$. Then $K$ contains a unique imaginary quadratic subfield $K_1$. Every imprimitive CM-type on $K$ is equivalent to $\Phi=G^{K_1}$. We have $\hat{K}=K_1$ and $\hat{\Phi}=\{\id\}$. There is a unique equivalence class of primitive CM-type. If $\tau$ is a generator of $\Gal(K/\bQ)$, a representative is $\Phi=\{\id,\tau,\tau^2\}$. We have $\hat{K}=K$ and $\hat{\Phi}=\{\id,\tau^{-1},\tau^{-2}\}$.

\case{v} $K$ is not Galois over $\bQ$ but $K=K^+(\sqrt{-D})$ for some square-free integer $D>0$. Then $G$ is a dihedral group of order $12$. There is an imprimitive equivalence class as in case \case{iv}, taking $K_1=\bQ(\sqrt{-D})$. There are primitive CM-types all of which are equivalent to the following one. Let $\tau$ be a generator of the unique cyclic subgroup of order $6$ of $G$, and put $\Phi=\{\id,\tau,\tau^2\}$. Then $\hat{K}=K$ and $\hat{\Phi}=\{\id, \tau^{-1},\tau^{-2}\}$.

In the remaining two cases, $K$ is not Galois over $\bQ$ and does not contain an imaginary quadratic subfield. Up to equivalence, there is a unique CM-type and it primitive. We always have $[\hat{K}:\bQ]=8$. 

\case{vi} $K^+$ is Galois over $\bQ$. Then $[L:\bQ]=24$ and $G^K$ is a Klein four-group. Furthermore, $G$ has four Sylow $3$-subgroups, and the restriction of the elements of any one of them to $K$ gives a CM-type $\Phi=\{\id,\tau,\tau^2\}$. Then $\hat{K}$ is the corresponding fixed field and $\hat{\Phi}$ is the restriction of the elements of $G^K$ to $\hat{K}$. 

\case{vii} $K^+$ is not Galois over $\bQ$. In this case, $[L:\bQ]=48$. Now $G^K$ is a dihedral group of order $8$ and we can take $\Phi$ to be the set of restrictions to $K$ of the elements $\{\id, \tau,\tau^2\}$ of one of the four Sylow $3$-subgroups of $G$. Then $G^{\hat{K}}$ is the unique subgroup of $G$ containing $\tau$ that is isomorphic to the symmetric group of degree $3$. Again $\hat{\Phi}$ is the set of restrictions to $\hat{K}$ of the elements of $G^K$.

Since there is only one primitive CM-type up to equivalence, we can test Estimate~\ref{Kfixedestimate} by determining all possible Weil numbers in the field $K$ without explicitly dealing with a reflex field, and then dividing out by the order of the automorphism group of $K$ (see Remark~\ref{autheuristics}). 

We wrote a program in Magma \cite{Magma} to compute the number $N(k,K,\rho_0,(a,b))$ of characteristic polynomials $C_\pi$ coming from triples $(r,\pi,p)$ with $a\leq r \leq b$ and $p\leq r^{\frac{\rho_0}{g}}$, and compared it with the value
\begin{equation*}
I(k,K,\rho_0,(a,b))=\frac{e(k,K)gw_Kh_{\hat{\Phi}}}{\#\big(\Aut(K)\big)\rho_0h_{\hat{K}}}\int_{a}^b\frac{du}{u^{2-\frac{\rho_0}{g}}(\log{u})^2}
\end{equation*} predicted by Estimate~\ref{Kfixedestimate}.  Looping over $r$, for each such prime pair $(r,p)$ satisfying $r\equiv 1 \pmod k$ and $p$ a primitive $k^{\text{th}}$ root of unity mod $r$, we search for $p$-Weil numbers in the following way:

\begin{enumerate}
\item 
Factorize $p\cO_K$ into prime ideals and make a list $D(p)$ of all possible ideal decompositions of the form $\ga\bar\ga = p\cO_K$ which are \emph{primitive}, that is, those for which there is no decomposition of the form $(\ga \cap K^0)(\bar\ga \cap K^0) = p\cO_{K^0}$ for any proper CM-subfield $K^0$.
\item 
For each pair $(\ga,\bar\ga) \in D(p)$ test whether $\ga$ is principal and if so find a generator $\gamma$.  Such an element satisfies $\gamma\bar\gamma = p\eta$ for some unit $\eta$ of $\cO_K$. 
Determine whether $\eta=\frac{\gamma\bar\gamma}{p}$ can be written in the form $\ve\bar\ve$.  If so then $\pi=\frac{\gamma}{\ve}$ is a $p$-Weil number and $\Gamma_\ga=\{\eta\pi \ :\ \eta \in U_K\}$ is the complete set of $p$-Weil numbers corresponding to $(\ga,\bar\ga)$.
\end{enumerate}
For each Weil number $\pi$ found, we check whether $r$ divides $\N_{K/\bQ}(\pi-1)$ and store the minimal polynomial $C_\pi$ (and its associated data $(r,p,\rho=g\frac{\log p}{\log r})$) for those $\pi$ satisfying this condition.

% When the degree and Galois group of the CM field is known in advance we are able to make some improvements to the above algorithm.  Most notably we can use congruence conditions modulo the conductor of the maximal abelian subfield of the Galois closure to avoid factoring ideals which could not possibly give rise to Weil numbers. 

Since $p$-Weil numbers are generators of principal ideals of the form $\N_{\hat{\Phi}}(\gP)$ where $\gP$ has norm $p$, we need only consider those $p$ for which there is a degree one prime above $p$ in $\hat{K}$.  We obtain necessary conditions by working in the maximal abelian subfield $M$ of the Galois closure $L$ of $K$.  We require that the decomposition field of a prime of $M$ above $p$ contains $F=\hat{K} \cap M$.  The Kronecker-Weber theorem tells us that $F$ is contained in $\bQ(\zeta_f)$ where $f$ is the conductor of $F$ and the decomposition group  $\Gal(\bQ(\zeta_f)/F) \triangleleft \Gal(\bQ(\zeta_f)/\bQ) \cong (\bZ/f\bZ)^*$  gives us congruence conditions on $p \pmod{f}$ for such a decomposition to occur.  For non-Galois CM-fields in  genus two or three, $F$ is a quadratic field and hence we need only compute ideal decompositions for half the congruence classes modulo $f$.  For Galois CM-fields of degree $2g\leq 6$ the Galois group is abelian and $\hat{K}=K=F$ so the proportion of primes we deal with is even smaller, namely $1/2g$.
In a similar manner, since we require that $r$ splits completely in $K$, we obtain further congruence restrictions on $r$ when the maximal abelian subextension of $K$ is not contained in $\bQ(\zeta_k)$.

We ran this program on a selection of quartic and sextic CM-fields for several values of $k$.  
The field invariants making up the constant term in the heuristic formula Estimate~\ref{Kfixedestimate} are varied in our sample.

For Galois CM-fields of degree $4$ and $6$ we computed the type norm map explicitly to determine the unique decomposition (up to the Galois action).  This approach of using the type norm map, while possible for other non-Galois fields, was slower than computing the ideal decompositions due to requiring to perform calculations in the larger Galois closure.

The tables have been placed together near the end of the paper for ease of use.  

Table \ref{table842} gives the values of $N(k,K,\rho_0,(10^4,5\cdot 10^5))$  with rho-values $\rho_0\leq 3.5$ and with several values of $k$, for the class number one CM-field $K=\bQ[X]/(X^4+4X^4+2)$ which is a cyclic Galois extension of $\bQ$.  Table \ref{cyclo5} presents a similar table for the field $\bQ(\zeta_5)$ which is another Galois cyclic CM-field of class number one. 
Table \ref{table12_8_13} presents a table in the same format, this time for a non-Galois quartic CM-field. The cyclic examples took between 5000-10000 seconds to compute, whereas the non-Galois example took 20000-30000 seconds.

Tables \ref{cyclo9}--\ref{Gal24:49_35_364_1183} give the values of $N(k,K,\rho_0,(10^4,5\cdot 10^5))$ with $\rho_0\leq 5.1$ for several values of $k$ for some sextic CM-fields belonging to cases \case{iv}, \case{v} and \case{vi} above.  It proved computationally too challenging to compute the heuristic value for a generic sextic field having Galois group of order $48$ so we did not produce any data for such a field.
The data for the Galois sextic CM-field examples was computed relatively quickly: approximately 2300 seconds for each value of $k$.  The sextic fields with $G$ of order $12$ took between 20000 and 60000 seconds each; the examples with $G$ of order $24$ took under 15000 seconds.  An explanation for why the order $24$ examples were quicker to compute than the order $12$ examples is that the latter have a proper CM-subfield (of degree $2$) so we must identify and discard the non-primitive ideal decompositions.

We would have liked to have extended the range of $r$, but the unavoidable large number of ideal factorizations in number fields prevented us from taking an interval for $r$ too large or high up. 

In almost all the cases we computed, there is good agreement between the computed value of $N(k,K,\rho_0,(a,b))$ and the expected value $I(k,K,\rho_0,(a,b))$. Noticeable exceptions occur in Table~\ref{table842}  when $k=2$ and $k=8$ and in Table~\ref{Gal12:229_24_144_27} when $k=24$, when the integral seems to seriously underestimate the actual number of triples found. To check whether the phenomenon persisted, we extended the computation to larger $r$ (up to $2\cdot 10^7$ in the cases $k=2$ and $k=8$ of Table~\ref{table842}) and found that the computed values were in much closer agreement with the expected ones.  

\subsection{An example in genus $4$}

We also computed pairing-friendly Weil polynomials for the non-Galois octic CM-field 
$$
K = \bQ[X]/(X^8 + 78X^6 + 1323X^4 + 7401X^2 + 9801) .
$$
Let $L$ be a Galois closure of $K$. Then $[L:\bQ]=24$ and $L$ contains a non-Galois sextic CM-field $K_6$ with Galois closure $L$.  Thus the Galois group $G=\Gal(L/\bQ)$ is that of case \case{vi} above.  It follows that the primitive CM-type on $K_6$ is a reflex CM-type for $K$.  Up to equivalence there are two primitive CM-types on $K$: $\Phi_6$ with reflex field $K_6$ and $\Phi_8$ with reflex field $K_8\cong K$.  There is also an imprimitive class of CM-types corresponding to extending the CM-type $\Phi_2$ of the imaginary quadratic field $K_2\cong\bQ(\zeta_3)$ contained in $K$.

Using the same method as described earlier we computed pairing-friendly primitive Weil polynomials of $K$.  We sorted them by CM-type $\Phi_i$ in order to compare the number $N_{\Phi_i}(k,K,\rho_0,(a,b))$ of pairing-friendly examples coming from $\Phi_i$ with the heuristic from Estimate~\ref{Kfixedestimate}:
\begin{equation*}
I_{\Phi_i}(k,K,\rho_0,(a,b))=\frac{e(k,K)gw_Kh_{\hat{\Phi}_i}}{\#\big(\Aut(K)\big)\rho_0h_{K_i}}\int_{a}^b\frac{du}{u^{2-\frac{\rho_0}{g}}(\log{u})^2}
\end{equation*} % predicted by Estimate~\ref{Kfixedestimate}.
(similar notations to before, now with $\Phi_i$ as a subscript).

Table \ref{tableg4} gives the values of $N_{\Phi_i}(k,K,\rho_0,(10^4,5\cdot 10^5))$ with $\rho_0 \leq 7.0$ for several values of $k$. 
%(all up the data for $k\in\{3,4,5,6\}$ took over $19$ days to compute).  
It turns out that 
$h_{\hat{\Phi}_6}/h_{K_6} = h_{\hat{\Phi}_8}/h_{K_8}$ in this example, so in fact $I_{\Phi_i}(k,K,\rho_0,(a,b))$ is the same for both primitive CM-types.

\begin{remark}
Since $w_K=w_{K_2}$ all imprimitive Weil numbers of $K$ are in $K_2$, and so an imprimitive Weil polynomial of $K$ with rho-value $\rho$ is a $4$th power of a Weil polynomial of $K_2$ with rho-value $\rho/4$.  As one would hope, the genus $1$ heuristic integral $I(k,K_2,\rho_0/4,(a,b))$ equals $I_{\Phi_2}(k,K,\rho_0,(a,b))$.
We confirmed that the imprimitive counts agree well with the heuristic estimate for $k=4,5$.  There were zero pairing-friendly examples found for the cases $k=3,6$.  These are two of the three ``exceptional'' cases, where the independence hypotheses are not satisfied, hence the heuristics do not apply.  See \cite{Bo} for details.
\end{remark}

\section{Asymptotics for a fixed maximal real subfield.} \label{fixedmaxreal}

When $\rho_0<g$, then the integral
\begin{equation*}
\int_2^{\infty}\frac{du}{u^{2-\frac{\rho_0}{g}}(\log{u})^2}
\end{equation*}
converges. Thus, the heuristic estimates Estimate~\ref{Kfixedestimate} suggest that, if $K$ is any CM-field of degree $2g$ and if $\rho_0< g$, there are only finitely many triples $(r,\pi,p)$ with rho-value less than $\rho_0$.

In order to try to understand where triples with rho-value less than $g$ might be located, we now develop a heuristic formula for the asymptotic growth of the number of triples with rho-value bounded above by $\rho_0$ and the CM-field $K$ varies but with a fixed maximal real subfield $K^+_0$. Thus, we fix $k$, $\rho_0$ and a totally real field $K^+_0$ of degree $g$ and seek an estimate for the number of triples $(r,\pi, p)$ with $r\leq x$, $p\leq r^{\frac{\rho_0}{g}}$ and $\pi$ is a $p$-Weil number lying in some CM-field whose maximal real subfield is $K^+_0$. We denote by $\cO^+_0$ the ring of integers of $K^+_0$. Furthermore, if $\alpha\in K^+_0$, we denote by $\{\alpha_i\mid 1\leq i\leq g\}$ the set of real embeddings of $\alpha$. 

Let $(r,\pi, p)$ be such a triple and write $\tau=\pi+p/\pi$. Then $\tau\in \cO^+_0$ and ($X$ denoting a variable)
\begin{equation*}
(X-\pi)(X-p/\pi)=X^2-\tau X+p\in \cO^+_0[X]\quad \text{ and }\quad \abs{\tau_i}\leq 2\sqrt{p} \text{ for all } i\in \{1,2,\dots, g\},
\end{equation*}   
the last inequalities being a consequence of the Weil bounds. Furthermore, the characteristic polynomial $C_{\pi}(X)$ of $\pi$ factors over $\bR[X]$ as
\begin{equation}  \label{Cpi}
C_{\pi}(X)=\prod_{i=1}^g(X^2-\tau_i X+p).
\end{equation}
Conversely, if $\tau\in \cO^+_0$ is such that $\abs{\tau_i}\leq 2\sqrt{p}$ for all $i$, then $X^2-\tau X+p$ has two roots $\pi$ and $p/\pi$ which are $p$-Weil numbers such that, if $\tau\neq \pm 2\sqrt{p}$, $K^+_0(\pi)$ is a CM field with maximal real subfield $K^+_0$. 

\begin{fixedmaxrealsubfieldformula} \label{Kplusfixedestimate}
Let $g\geq 1$ be an integer, let $K^+_0$ be a totally real field of degree $g$, let $k\geq 2$ be an integer and let $\rho_0> \max(1,\frac{g}{\varphi(k)})$ be a real number with $\rho_0\neq \frac{2g}{g+2}$. Then the number $R(k,K^+_0,\rho_0,x)$ of triples $(r,\pi,p)$ with $r\leq x$ is asymptotically equivalent as $x\to \infty$ to
\begin{equation*}
R(k,K^+_0,\rho_0,x)\sim \frac{g4^{g+1}e(k,K^+_0)}{\rho_0(g+2)d_0^{1/2}}\int_{2}^x\frac{u^{\rho_0\big(\frac{1}{2}+\frac{1}{g}\big)-2}du}{(\log{u})^2}.
\end{equation*} 
Here $d_0$ denotes the discriminant of $K^+_0$ and $e(k,K^+_0)$ the degree of $K_0^+\cap \bQ(\zeta_k)$ over $\bQ$. 
\end{fixedmaxrealsubfieldformula}

Here the integral diverges to $\infty$ with $x$ if and only if $\rho_0\big(\frac{1}{2}+\frac{1}{g}\big)> 1$. In particular, if $g=1$ or if $g=2$, we expect $R(k,K^+_0,\rho_0,x)$ to tend to $\infty$ with $x$ for all $\rho_0>1$ but, as $g$ increases, we expect this to be true only for $\rho_0$ larger than a bound which approaches $2$.

\begin{remark}\label{kplusauthheuristics}
Let $(r,\pi, p)$ be a triple as above, and let $\tau=\pi+p/\pi$. If $\pi$ is real, then $\pi=\pm \sqrt{p}$ and $\pi$ must belong to $K_0^+$. This can occur for only finitely many $p$. Otherwise, $K_0^+(\pi)$ is uniquely determined by $\pi$, and $\tau$ arises from the two $p$-Weil numbers $\pi$ and $p/\pi$. It follows that the number of triples $(r,\pi,p)$ that correspond to the same triple $(r,C,p)$ with $C$ a characteristic polynomial is equal to twice the number of conjugates of $\tau$; in particular if $K_0^+=\bQ(\tau)$, it is equal to $2\# \big(\Aut(K_0^+)\big)$. The triples are ordinary when $p$ is unramified in $K_0^+(\pi)$.
\end{remark}

We now indicate a heuristic argument that leads to Estimate~\ref{Kplusfixedestimate}. From elementary results about the geometry of algebraic number fields, we know that as $T\to \infty$, the number $R(K^+_0,T)$ of $\tau\in \cO^+_0$ such that $\abs{\tau_i}\leq T$ for all $i$ satisfies
\begin{equation*}
R_0(K^+_0,T)\sim (2T)^gd_0^{-1/2}.
\end{equation*}

As before, the probability that $r$ divides $\Phi_k(p)$ with $r$ prime is $\frac{1}{r\log{r}}$. On the other hand, by the prime ideal theorem in number fields the expected number of degree one prime ideals of $K_0^+$ dividing $r$ given that $r$ splits in $\bQ(\zeta_k)$ is equal to $e(k,K^+_0)$. In view of this, we assume that the probability that an integer $r$ is prime and divides both $\Phi_k(p)$ and $\N_{K^+_0(\pi)/\bQ}(\pi-1)$ is equal to $\frac{e(k,K_0^+)}{r^2\log{r}}$. 

Thus we expect the number $R(k,K^+_0,\rho_0,x)$ of triples $(r,\pi,p)$ with $r\leq x$ and $p\leq r^{\frac{\rho_0}{g}}$ to be equivalent to
\begin{equation}
\sum_{r\leq x}\frac{e(k,K_0^+)}{r^2\log{r}}\sum_{p\leq r^{\frac{\rho_0}{g}}}2R_0(K^+_0,2\sqrt{p}),
\end{equation}
where the sum over $r$ is over integers and that over $p$ is over primes, and the $2$ appears before the $R_0(K^+_0,2\sqrt{p})$ because we distinguish between $\pi$ and $p/\pi$. Hence
\begin{equation}\label{eqN1}
R(k,K^+_0,\rho_0,x)\sim \frac{2\cdot 4^ge(k,K^+_0)}{d_0^{1/2}}\sum_{r\leq x}\frac{1}{r^2\log{r}}\sum_{p\leq r^{\frac{\rho_0}{g}}}p^{\frac{g}{2}}
\end{equation} 

Now, it follows from the prime number theorem by an easy argument using Abel summation that, if $\alpha\geq 0$, then the sum over primes $\sum_{p\leq U}p^\alpha$ is asymptotically equivalent to $\frac{U^{\alpha+1}}{(\alpha+1)\log{U}}$ as $U\to \infty$ (see for example \cite{Ld}, pages 203-205 for a more general statement). Applying this with $\alpha=\frac{g}{2}$ and $U=r^{\frac{\rho_0}{g}}$, (\ref{eqN1}) becomes 
\begin{equation}
R(k,K^+_0,\rho_0,x)\sim \frac{g4^{g+1}e(k,K^+_0)}{\rho_0(g+2)d_0^{1/2}}\sum_{r\leq x}\frac{r^{\frac{\rho_0}{g}\big(\frac{g}{2}+1\big)-2}}{(\log{r})^2}.
\end{equation}

Replacing the sum by an integral and rearranging slightly leads to Estimate~\ref{Kplusfixedestimate}. 

\begin{remark} \label{clusterrmk}
Suppose that there are $C$ prime ideals of norm $r$ in $K_0^+$ (where $0\leq C\leq g$). Then the number of $\tau\in \cO^+_0$ with $\abs{\tau}_i\leq 2\sqrt{p}$ for all embeddings $i$ of $\tau$ in $\bR$ is equivalent to
\begin{equation*}
4^gp^{\frac{g}{2}}d_0^{-1/2}.
\end{equation*} 
On the other hand, if $\gr^+$ is a prime ideal of norm $r$ in $K_0^+$, the probability that $\tau\equiv p+1\pmod{\gr^+}$ is $\frac{1}{r}$. Hence the number of elements $\tau$ in this range and $\gr^+$ a degree one prime ideal of $K_0^+$ dividing $r$ such that $\tau\equiv p+1\pmod{\gr^+}$ should be roughly
\begin{equation*}
4^gp^{\frac{g}{2}}Cr^{-1}d_0^{-1/2}.
\end{equation*}     
In particular, if $p$ is close to $r^{\frac{\rho_0}{g}}$, this is close to
\begin{equation} \label{clusterrmkform}
4^gr^{\frac{\rho_0}{2}-1}Cd_0^{-1/2}.
\end{equation}

Thus, if $\rho_0>2$, we expect that when $r$ is large, and $p$ is a prime close to $r^{\frac{\rho_0}{g}}$, the number of $p$-Weil numbers $\pi$ with $\pi+p/\pi\in K^+_0$ is close to (\ref{clusterrmkform}), which tends to infinity with $r$.

For a numerical illustration, see Remark \ref{clusterexample}. 
\end{remark}

\section{Numerical evidence in the fixed maximal real subfield case}  \label{fixedmaxrealcomp}

We continue to use the notation introduced in the previous section.  For small $x$ and $\rho_0$, we can compute $R(k,K^+_0,\rho_0,x)$ as follows. Since $r\leq x$, we know that $\abs{\tau_i}\leq 2\sqrt{p}\leq 2r^{\frac{\rho_0}{2g}}\leq 2x^{\frac{\rho_0}{2g}}$. Hence, we need to:

\begin{enumerate}
\item[(1)] Make a list $\cL$ of all $\tau\in \cO^+_0$ such that $\abs{\tau_i}\leq 2x^{\frac{\rho_0}{2g}}$ for all $i$;

\item[(2)] For each $\tau\in \cL$, factor $\Phi_k(\tau-1)$ into prime ideals in $K^+_0$ and make a list $\cM(\tau)$ of all degree one primes $\gr^+$ dividing $\Phi_k(\tau-1)$ of norm $r$ such that $x\geq r\geq \big(\frac{\abs{\tau_i}}{2}\big)^{\frac{2g}{\rho_0}}$ for all $i$;

\item[(3)] For each $\gr^+\in M(\tau)$, search for primes $p\leq x^{\frac{\rho_0}{g}}$ such that $p\equiv \tau-1\pmod{\gr^+}$ and $|\tau_i|\leq 2\sqrt{p}$ for all $i$. 
\end{enumerate}

The condition $p\equiv \tau-1\pmod{\gr^+}$ of (3) ensures that $r$ divides $\Phi_k(p)$. Thus, for any $\tau\in \cL$, $\gr^+\in M(\tau)$ and prime $p$ as in (3), the triples $(r,\pi,p)$ and $(r,p/\pi,p)$ with $\pi$ and $p/\pi$ the roots of $X^2-\tau X+p$ contribute towards the total $R(k,K^+_0,\rho_0,x)$. 

Conversely, if $(r,\pi,p)$ is a triple with $r\leq x$ and $p\leq r^{\frac{\rho_0}{g}}$, and if $\tau=\pi+p/\pi$, then $\tau\in \cL$. Since $r$ divides $\N_{K/\bQ}(\pi-1)$, some prime ideal $\gr^+$ of $K^+$ above $r$ must divide $p+1-\tau$. Then $\gr^+\in M(\tau)$ and $p$ satisfies (3), so that $(r,\pi,p)$ will be detected by the above search. 

Of course the major drawback of this approach is the need to factor $\Phi_k(\tau-1)$. Since the size of $\Phi_k(\tau-1)$ depends on the degree of $\Phi_k$, it is necessary to choose $k$ with $\varphi(k)$ small. On the other hand, decreasing $\rho$ reduces the size of the list $\cL$ of step (1).  In any case, in practice it is only possible to make meaningful computations when $\varphi(k)$ and $x$ are small. 

Using this method, we computed a few tables with $\rho_0$ at most equal to $g$.

Let $R_c(k,K^+_0,\rho_0,(a,b))$ be the number of distinct irreducible characteristic polynomials (\ref{Cpi})  associated to Weil numbers $\pi$ belonging to triples $(r,\pi,p)$ with $a\leq r\leq b$, $\pi+p/\pi\in K_0^+$ and $p\leq r^{\frac{\rho_0}{g}}$. Tables \ref{Real_sqrt2}--\ref{real_g3_cyclo7+} compare the values of $R_c(k,K^+_0,\rho_0,(a,b))$ with the heuristic estimate
\begin{equation*}
J=J(k,K^+_0,\rho_0,(a,b))=\frac{g4^{g+1}e(k,K^+_0)}{2\#\big(\Aut(K_0^+)\big)\rho_0(g+2)d_0^{1/2}}\int_{a}^b\frac{u^{\rho_0\big(\frac{1}{2}+\frac{1}{g}\big)-2}du}{(\log{u})^2}.
\end{equation*}

Table \ref{Real_sqrt2} shows the values  of $R_c(k,\bQ(\sqrt{2}),\rho_0,(10^3,10^5))$ for $k\in \{3,4,5,6,7,8,12\}$ and values of $\rho_0$ between $1$ and $2$. On the other hand, Table \ref{sqrts2.0} presents the values of $R_c(k,\bQ(\sqrt{d}),2.0,(10^3,10^5))$ for $k \in \{3,4,5,6,12\}$ and for squarefree $d\leq 50$.   The entries for which $e(k,\bQ(\sqrt{d}))=2$ are marked with an asterisk; the predicted value for these entries is $2J$. In all other cases, $e(k,\bQ(\sqrt{d}))=1$ and the predicted value is $J$. Finally, Table \ref{real_g3_cyclo7+} shows the values  of $R_c(k,\bQ(\zeta_7+\zeta_7^{-1}),\rho_0,(10^3,10^4))$ for $k\in \{3,4,5,6,7\}$ and values of $\rho_0$ between $1.5$ and $3$.

The running times in Table \ref{sqrts2.0} were dependent on the size of $\varphi(k)$, the degree of the $k^{\text{th}}$ cyclotomic polynomial.  For $k\neq 5,12$ each table entry took under 10 minutes to compute.  The value which took the most time to compute was $k=5$ for the field $\bQ(\sqrt{5})$ which just under $3\frac{1}{2}$ hours. The computations in Table \ref{real_g3_cyclo7+} took under three hours when $k\in \{3,4,6\}$ but over ten times longer when $k=5$ and $k=7$.

\begin{remark}
In almost all cases in Tables~\ref{Real_sqrt2}, \ref{sqrts2.0} and \ref{real_g3_cyclo7+}, the CM-field $K_0^+(\pi)$ has Galois closure of maximal size. For instance, over $10^4$ examples are listed in Table~\ref{Real_sqrt2}, but in only $58$ of them is the field $\bQ(\sqrt{2},\pi)$ Galois; in $54$ cases the field is biquadratic and in the other four it is cyclic.
\end{remark}

\begin{remark} \label{clusterexample}
We can also use Table  \ref{real_g3_cyclo7+} to illustrate Remark \ref{clusterrmk}. When $K_0^+=\bQ(\zeta_7+\zeta_7^{-1})$ and $k=5$, the table shows that there are no triples with rho-value between $2.3$ and $2.4$, and $183$ with rho-value between $2.4$ and $2.5$. However, only three values of $(r,p)$ account for all these triples: $(1051,307)$, $(5741,1229)$ and $(6091,1321)$ which give rise respectively to $46$, $66$ and $74$ triples with corresponding rho-values $2.469$, $2.466$ and $2.474$. If we substitute $r=1051$, $5741$ or $6091$ in  (\ref{clusterrmkform}) with $g=3$, $C=3$ and divide by $\#\Aut(\bQ(\zeta_7+\zeta_7^{-1}))=3$, we find that the predicted number of triples is respectively $46.9$, $68.7$ and $72.1$.  
\end{remark}

\include{tables-plain}

\end{document}

%% file: tables-plain.tex
\hyphenation{Invariants}
\newcolumntype{H}{>{\setbox0=\hbox\bgroup}c<{\egroup}@{}}
\newcolumntype{C}{>{$}c<{$}}
\newcolumntype{T}{@{\hspace{0.23\tabcolsep}}>{$}c<{$}@{\hspace{0.23\tabcolsep}}}
\newcolumntype{D}{@{\hspace{0.8\tabcolsep}}>{$}c<{$}@{\hspace{0.8\tabcolsep}}}

\newcommand\blockcomment[1]{}

\begin{sidewaystable}[p]
\vskip12cm

\begin{center}
\begin{tabular}{|D||DDDDDDD|DDD|DDD|}
\hline
\rho_0 & k=2 & k=3 & k=4 & k=5 & k=6 & k=7 & I & k=8 & k=24 & I & k=16 & k=32 & I \\ \hline
2.8 & 2 & 3 & 1 & 0 & 0 & 0 & 1.02 & 7 & 1 & 2.03 & 3 & 4 & 4.07 \\ \hline
2.9 & 4 & 3 & 2 & 0 & 3 & 1 & 1.74 & 8 & 1 & 3.48 & 7 & 5 & 6.97 \\ \hline
3.0 & 8 & 3 & 6 & 1 & 5 & 2 & 3.00 & 16 & 3 & 6.00 & 10 & 11 & 11.99 \\ \hline
3.1 & 14 & 5 & 8 & 2 & 10 & 3 & 5.18 & 20 & 5 & 10.36 & 22 & 17 & 20.73 \\ \hline
3.2 & 22 & 9 & 9 & 6 & 13 & 5 & 8.99 & 23 & 15 & 17.98 & 43 & 33 & 35.96 \\ \hline
3.3 & 30 & 14 & 15 & 12 & 26 & 14 & 15.66 & 36 & 30 & 31.31 & 63 & 58 & 62.62 \\ \hline
3.4 & 46 & 27 & 26 & 23 & 40 & 31 & 27.37 & 61 & 55 & 54.73 & 112 & 104 & 109.46 \\ \hline
3.5 & 68 & 51 & 59 & 38 & 59 & 49 & 48.00 & 99 & 110 & 96.00 & 178 & 187 & 192.00 \\ \hline
\end{tabular}
\caption{ 
Values of $N(k,K,\rho_0,(10^4,5\cdot 10^5))$ for $K=\bQ[X]/(X^4+4X^2+2)$. Invariants: $w_K=2$, $h_{\hat\Phi}=h_{\hat K}=1$, $G$ cyclic. 
}
\label{table842}

\begin{tabular}{|D||DDDDDDD|DDDDDD|}
\hline
\rho_0 &   k=2 &   k=3 &    k=4 & k=12 & k=24 & k=36 &  I & k=5 &  k=10 & k=15 & k=20 & k=25 &  I \\ \hline
2.5 & 0 & 3 & 0 & 2 & 2 & 2 & 1.04 & 2 & 4 & 9 & 2 & 4 & 4.15 \\ \hline
2.6 & 2 & 3 & 2 & 3 & 2 & 6 & 1.75 & 6 & 10 & 12 & 3 & 6 & 7.01 \\ \hline
2.7 & 2 & 5 & 2 & 3 & 4 & 7 & 2.98 & 10 & 22 & 17 & 5 & 6 & 11.91 \\ \hline
2.8 & 2 & 6 & 2 & 6 & 6 & 10 & 5.08 & 14 & 26 & 29 & 14 & 9 & 20.33 \\ \hline
2.9 & 6 & 9 & 8 & 8 & 9 & 10 & 8.71 & 26 & 46 & 45 & 32 & 22 & 34.84 \\ \hline
3.0 & 10 & 15 & 14 & 18 & 17 & 18 & 14.99 & 64 & 70 & 72 & 49 & 51 & 59.97 \\ \hline
3.1 & 16 & 27 & 20 & 32 & 24 & 27 & 25.91 & 106 & 124 & 125 & 83 & 93 & 103.63 \\ \hline
3.2 & 26 & 44 & 43 & 52 & 35 & 50 & 44.95 & 176 & 168 & 210 & 150 & 162 & 179.79 \\ \hline
3.3 & 70 & 76 & 72 & 82 & 72 & 87 & 78.28 & 302 & 302 & 335 & 282 & 319 & 313.12 \\ \hline
3.4 & 112 & 142 & 140 & 143 & 130 & 141 & 136.83 & 574 & 560 & 597 & 534 & 578 & 547.30 \\ \hline
3.5 & 212 & 250 & 241 & 258 & 235 & 251 & 240.00 & 1000 & 1000 & 1049 & 977 & 1006 & 959.99 \\ \hline
\end{tabular}
\caption{
Values of $N(k,K,\rho_0,(10^4,5\cdot 10^5))$ for $K=\mathbb{Q}(\zeta_5)$. Invariants:  $w_K=10$, $h_{\hat\Phi}=h_{\hat K}=1$, $G$ cyclic. 
}
\label{cyclo5}
\begin{tabular}{|D||DDDDDDD|DDDD|}
\hline
\rho_0 &  k=2 &   k=3 &  k=4 & k=5 &  k=6 &  k=7 &  I &  k=12 & k=24 & k=36 &  I \\ \hline
2.7 & 2 & 0 & 0 & 1 & 2 & 2 & 1.19 & 2 & 2 & 2 & 2.38 \\ \hline
2.8 & 2 & 2 & 2 & 4 & 3 & 3 & 2.03 & 2 & 4 & 6 & 4.07 \\ \hline
2.9 & 6 & 5 & 3 & 6 & 3 & 4 & 3.48 & 8 & 8 & 9 & 6.97 \\ \hline
3.0 & 6 & 8 & 6 & 10 & 6 & 7 & 6.00 & 17 & 14 & 11 & 11.99 \\ \hline
3.1 & 8 & 13 & 11 & 11 & 10 & 14 & 10.36 & 25 & 25 & 17 & 20.73 \\ \hline
3.2 & 16 & 23 & 19 & 20 & 17 & 25 & 17.98 & 44 & 43 & 36 & 35.96 \\ \hline
3.3 & 32 & 31 & 26 & 34 & 27 & 39 & 31.31 & 65 & 71 & 64 & 62.62 \\ \hline
3.4 & 58 & 59 & 56 & 57 & 54 & 66 & 54.73 & 116 & 116 & 115 & 109.46 \\ \hline
3.5 & 100 & 97 & 93 & 93 & 96 & 117 & 96.00 & 206 & 195 & 191 & 192.00 \\ \hline
\end{tabular}
\caption{
Values of $N(k,K,\rho_0,(10^4,5\cdot 10^5))$ for  $K=\bQ[X]/(X^4+8X^2+13)$. Invariants:  $w_K=2$, $h_{\hat\Phi}=h_{\hat K}=2$, $G$ dihedral. 
}
\label{table12_8_13}
\end{center}
\end{sidewaystable}

\begin{table}[p]
\begin{tabular}{|D||DDDD|DDD|DDD|}
\hline
\rho_0 & k=2 & k=4 & k=5 & I & k=3 & k=6 & I & k=9 & k=18 & I \\ \hline
4.0 & 6 & 3 & 0 & 2.99 & 2 & 4 & 5.99 & 22 & 18 & 17.97 \\ \hline
4.1 & 8 & 6 & 2 & 4.27 & 6 & 8 & 8.54 & 34 & 24 & 25.62 \\ \hline
4.2 & 10 & 6 & 6 & 6.10 & 10 & 18 & 12.20 & 46 & 44 & 36.60 \\ \hline
4.3 & 14 & 10 & 8 & 8.73 & 14 & 22 & 17.46 & 64 & 54 & 52.38 \\ \hline
4.4 & 16 & 11 & 13 & 12.52 & 20 & 30 & 25.04 & 82 & 72 & 75.13 \\ \hline
4.5 & 24 & 15 & 23 & 17.99 & 30 & 38 & 35.98 & 124 & 116 & 107.94 \\ \hline
4.6 & 32 & 24 & 30 & 25.90 & 50 & 62 & 51.79 & 180 & 160 & 155.37 \\ \hline
4.7 & 44 & 34 & 42 & 37.34 & 80 & 80 & 74.68 & 260 & 236 & 224.05 \\ \hline
4.8 & 68 & 51 & 62 & 53.94 & 114 & 116 & 107.88 & 390 & 330 & 323.63 \\ \hline
4.9 & 90 & 71 & 82 & 78.04 & 166 & 162 & 156.09 & 568 & 454 & 468.27 \\ \hline
5.0 & 136 & 104 & 114 & 113.11 & 250 & 224 & 226.22 & 812 & 658 & 678.66 \\ \hline
5.1 & 224 & 169 & 159 & 164.19 & 380 & 328 & 328.38 & 1238 & 944 & 985.15 \\ \hline
\end{tabular}
\caption{
Values of $N(k,K,\rho_0,(10^4,5\cdot 10^5))$ for $K=\bQ(\zeta_9)$. Invariants:  $w_K=18$, $h_{\hat\Phi}=h_{\hat K}=1$, $G$ cyclic. 
}
\label{cyclo9}

\begin{tabular}{|D||DDDDD|DDDD|}
\hline
\rho_0 &  k=2 &  k=4 &  k=5 & k=32 &  I & k=3 & k=6 & k=24 &   I \\ \hline
3.9 & 0 & 3 & 0 & 0 & 1.05 & 2 & 4 & 3 & 2.10 \\ \hline
4.0 & 0 & 3 & 0 & 0 & 1.50 & 2 & 4 & 5 & 2.99 \\ \hline
4.1 & 0 & 3 & 0 & 1 & 2.13 & 4 & 6 & 7 & 4.27 \\ \hline
4.2 & 2 & 3 & 0 & 2 & 3.05 & 6 & 6 & 10 & 6.10 \\ \hline
4.3 & 4 & 5 & 0 & 4 & 4.37 & 8 & 6 & 15 & 8.73 \\ \hline
4.4 & 6 & 5 & 2 & 6 & 6.26 & 14 & 8 & 21 & 12.52 \\ \hline
4.5 & 12 & 8 & 6 & 9 & 9.00 & 20 & 14 & 32 & 17.99 \\ \hline
4.6 & 16 & 12 & 9 & 13 & 12.95 & 22 & 24 & 53 & 25.90 \\ \hline
4.7 & 22 & 15 & 13 & 20 & 18.67 & 32 & 34 & 67 & 37.34 \\ \hline
4.8 & 40 & 23 & 24 & 30 & 26.97 & 44 & 50 & 84 & 53.94 \\ \hline
4.9 & 50 & 35 & 32 & 42 & 39.02 & 62 & 80 & 119 & 78.04 \\ \hline
5.0 & 64 & 52 & 57 & 58 & 56.55 & 110 & 118 & 160 & 113.11 \\ \hline
5.1 & 88 & 74 & 96 & 84 & 82.10 & 164 & 170 & 214 & 164.19 \\ \hline
\end{tabular}
\caption{
Values of $N(k,K,\rho_0,(10^4,5\cdot 10^5))$ for $K=\bQ[X]/(X^6+24X^4+144X^2+27)$. Invariants:  $w_K=6$, $h_{\hat\Phi}=1$, $h_{\hat K}=2$, $G$ of order $12$. 
}
\label{Gal12:229_24_144_27}

\begin{tabular}{|D||DDDDDD|DDDD|}
\hline
\rho_0 &  k=2 &  k=3 &  k=4 &  k=5 &  k=6 &  I & k=7 & k=14 & k=35 & I \\ \hline
4.4 & 0 & 2 & 0 & 2 & 3 & 1.04 & 5 & 2 & 4 & 3.13 \\ \hline
4.5 & 0 & 2 & 0 & 2 & 4 & 1.50 & 10 & 4 & 4 & 4.50 \\ \hline
4.6 & 2 & 2 & 0 & 3 & 5 & 2.16 & 11 & 5 & 6 & 6.47 \\ \hline
4.7 & 2 & 3 & 0 & 4 & 6 & 3.11 & 15 & 7 & 10 & 9.34 \\ \hline
4.8 & 2 & 6 & 3 & 6 & 8 & 4.49 & 16 & 14 & 11 & 13.48 \\ \hline
4.9 & 2 & 8 & 4 & 8 & 8 & 6.50 & 23 & 23 & 17 & 19.51 \\ \hline
5.0 & 8 & 13 & 6 & 15 & 10 & 9.43 & 37 & 37 & 25 & 28.28 \\ \hline
5.1 & 12 & 14 & 9 & 18 & 14 & 13.68 & 48 & 49 & 40 & 41.05 \\ \hline
\end{tabular}
\caption{
Values of $N(k,K,\rho_0,(10^4,5\cdot 10^5))$ for $K=\bQ[X]/(X^6+35X^4+364X^2+1183)$. Invariants:  $w_K=2$, $h_{\hat\Phi}=4$, $h_{\hat K}=16$, $G$ of order $24$. 
}
\label{Gal24:49_35_364_1183}

\end{table}

\newcolumntype{B}{>{$}c<{$}}
\begin{table}[h]
\centering
\begin{tabular}{|B||BBBBB|BBBBB|}
\hline
 & \multicolumn{2}{B}{k=4} & \multicolumn{2}{B}{k=5} & \text{heuristic} & \multicolumn{2}{B}{k=3} & \multicolumn{2}{B}{k=6} & \text{heuristic} \\ 
\rho_0 & N_{\Phi_6} & N_{\Phi_8} & N_{\Phi_6} & N_{\Phi_8} & I_{\Phi_6}\!=\!I_{\Phi_8} & N_{\Phi_6} & N_{\Phi_8} & N_{\Phi_6} & N_{\Phi_8} & I_{\Phi_6}\!=\!I_{\Phi_8}\\ \hline
6.0 &   5 &   9 &  16 &  12 &   9.00 &  16 &  20 &  18 &  14 &  18.00 \\ \hline
6.1 &   6 &  11 &  18 &  19 &  11.82 &  20 &  24 &  26 &  20 &  23.64 \\ \hline
6.2 &  12 &  14 &  21 &  26 &  15.54 &  30 &  28 &  36 &  28 &  31.09 \\ \hline
6.3 &  21 &  25 &  27 &  32 &  20.47 &  42 &  38 &  56 &  38 &  40.93 \\ \hline
6.4 &  31 &  39 &  32 &  37 &  26.97 &  56 &  62 &  74 &  50 &  53.94 \\ \hline
6.5 &  40 &  51 &  41 &  46 &  35.57 &  68 &  74 &  94 &  62 &  71.15 \\ \hline
6.6 &  49 &  64 &  53 &  55 &  46.96 &  90 &  96 & 128 &  82 &  93.94 \\ \hline
6.7 &  62 &  81 &  74 &  72 &  62.07 & 136 & 130 & 152 & 116 & 124.14 \\ \hline
6.8 &  85 & 104 &  89 &  94 &  82.10 & 176 & 176 & 196 & 152 & 164.19 \\ \hline
6.9 & 117 & 133 & 118 & 131 & 108.68 & 240 & 216 & 236 & 222 & 217.36 \\ \hline
7.0 & 157 & 167 & 159 & 171 & 144.00 & 300 & 286 & 300 & 314 & 288.00 \\ \hline
\end{tabular}
\caption{ 
Values of $N_{\Phi_i}(k,K,\rho_0,(10^4,5\cdot 10^5))$
for the field $K=\bQ[X]/(X^8 + 78X^6 + 1323X^4 + 7401X^2 + 9801)$.
Invariants: $w_K=6$, $h_{\hat{\Phi}_6}=4, h_{\hat K_6}=8$, $h_{\hat{\Phi}_8}=2, h_{\hat K_8}=4$. 
}
\label{tableg4}

\begin{tabular}{|D||DDDDDDD|DD|}
\hline
\rho_0 &  k=3 &  k=4 &  k=5 &  k=6 &  k=7 & k=12 &  J & k=8 &  J \\ \hline
1.0 & 1 & 0 & 0 & 0 & 0 & 1 & 0.16 & 0 & 0.33 \\ \hline
1.1 & 1 & 0 & 0 & 1 & 0 & 1 & 0.36 & 0 & 0.73 \\ \hline
1.2 & 2 & 0 & 0 & 1 & 2 & 2 & 0.83 & 1 & 1.65 \\ \hline
1.3 & 4 & 1 & 0 & 1 & 3 & 3 & 1.92 & 1 & 3.85 \\ \hline
1.4 & 7 & 2 & 5 & 5 & 4 & 6 & 4.59 & 7 & 9.18 \\ \hline
1.5 & 15 & 11 & 14 & 15 & 12 & 17 & 11.21 & 22 & 22.42 \\ \hline
1.6 & 36 & 22 & 28 & 34 & 25 & 37 & 27.95 & 62 & 55.90 \\ \hline
1.7 & 81 & 68 & 62 & 88 & 62 & 80 & 71.04 & 157 & 142.09 \\ \hline
1.8 & 200 & 194 & 192 & 219 & 161 & 210 & 183.80 & 384 & 367.60 \\ \hline
1.9 & 493 & 518 & 467 & 496 & 534 & 543 & 483.16 & 940 & 966.33 \\ \hline
2.0 & 1346 & 1418 & 1267 & 1331 & 1295 & 1321 & 1288.45 & 2572 & 2576.91 \\ \hline
\end{tabular}
\caption{Values of $R_c(k,K^+_0,\rho_0,(10^3,10^5))$ for $K^+_0 = \bQ(\sqrt{2})$.}
\label{Real_sqrt2}

\begin{tabular}{|T||TTTTTT||T||TTTTTT|}
\hline
d & k\!=\!3 & k\!=\!4 & k\!=\!5 & k\!=\!6 & k\!=\!12 & J & 
d & k\!=\!3 & k\!=\!4 & k\!=\!5 & k\!=\!6 & k\!=\!12 & J \\ \hline
2 & 1346 & 1418 & 1267 & 1331 & 1321 & 1288.45 & 26 & 365 & 408 & 368 & 374 & 358 & 357.35 \\ \hline
3 & 1144 & 1093 & 1049 & 1103 & 2199^* & 1052.02 & 29 & 675 & 718 & 688 & 662 & 660 & 676.73 \\ \hline
5 & 1650 & 1808 & 3306^* & 1670 & 1703 & 1629.78 & 30 & 356 & 338 & 322 & 346 & 354 & 332.68 \\ \hline
6 & 789 & 794 & 774 & 753 & 751 & 743.89 & 31 & 351 & 351 & 333 & 345 & 328 & 327.27 \\ \hline
7 & 755 & 718 & 634 & 667 & 708 & 688.71 & 33 & 643 & 687 & 621 & 664 & 640 & 634.39 \\ \hline
10 & 659 & 635 & 573 & 599 & 616 & 576.21 & 34 & 325 & 324 & 336 & 287 & 291 & 312.50 \\ \hline
11 & 574 & 580 & 534 & 553 & 567 & 549.40 & 35 & 319 & 341 & 285 & 311 & 349 & 308.00 \\ \hline
13 & 1090 & 1043 & 1064 & 975 & 1084 & 1010.75 & 37 & 634 & 596 & 654 & 614 & 609 & 599.12 \\ \hline
14 & 521 & 526 & 494 & 491 & 432 & 486.99 & 38 & 309 & 320 & 299 & 313 & 302 & 295.59 \\ \hline
15 & 486 & 460 & 487 & 443 & 475 & 470.48 & 39 & 325 & 334 & 280 & 307 & 306 & 291.78 \\ \hline
17 & 967 & 954 & 952 & 880 & 902 & 883.87 & 41 & 609 & 651 & 580 & 537 & 602 & 569.14 \\ \hline
19 & 422 & 480 & 450 & 395 & 412 & 418.03 & 42 & 320 & 280 & 316 & 303 & 255 & 281.16 \\ \hline
21 & 883 & 753 & 799 & 798 & 810 & 795.25 & 43 & 302 & 300 & 296 & 274 & 300 & 277.88 \\ \hline
22 & 396 & 415 & 405 & 379 & 414 & 388.48 & 46 & 307 & 289 & 258 & 300 & 253 & 268.66 \\ \hline
23 & 377 & 393 & 418 & 378 & 396 & 379.94 & 47 & 273 & 258 & 311 & 257 & 252 & 265.79 \\ \hline
\end{tabular}
\caption{
Values of $R_c(k,\bQ(\sqrt{d}),2.0,(10^3,10^5))$ for $k\!\in\!\{3,4,5,6,12\}$ and $d\!\leq\! 50$ squarefree. Asterisks indicate the cases where $e(k,\bQ(\sqrt{d})=2$.}
\label{sqrts2.0}

\end{table}

\begin{table}[p]

\begin{tabular}{|D||DDDDD|DD|}
\hline
\rho_0 & k=3 & k=4 & k=5 & k=6 & J & k=7 & J \\ \hline
1.5 & 3 & 0 & 1 & 0 & 0.65 & 2 & 1.96 \\ \hline
1.6 & 3 & 0 & 1 & 1 & 1.20 & 2 & 3.60 \\ \hline
1.7 & 10 & 11 & 1 & 3 & 2.22 & 6 & 6.66 \\ \hline
1.8 & 10 & 11 & 1 & 5 & 4.14 & 9 & 12.41 \\ \hline
1.9 & 10 & 28 & 1 & 9 & 7.75 & 24 & 23.26 \\ \hline
2.0 & 18 & 42 & 1 & 15 & 14.61 & 30 & 43.84 \\ \hline
2.1 & 32 & 53 & 12 & 35 & 27.70 & 77 & 83.10 \\ \hline
2.2 & 144 & 82 & 40 & 68 & 52.78 & 230 & 158.33 \\ \hline
2.3 & 197 & 82 & 97 & 160 & 101.05 & 324 & 303.15 \\ \hline
2.4 & 244 & 232 & 97 & 236 & 194.37 & 716 & 583.11 \\ \hline
2.5 & 354 & 519 & 280 & 362 & 375.53 & 1028 & 1126.60 \\ \hline
2.6 & 557 & 1048 & 714 & 865 & 728.59 & 1647 & 2185.76 \\ \hline
2.7 & 1211 & 1654 & 1314 & 1132 & 1419.19 & 3267 & 4257.58 \\ \hline
2.8 & 2474 & 3050 & 2640 & 1598 & 2774.87 & 9820 & 8324.62 \\ \hline
2.9 & 5136 & 5527 & 5330 & 3993 & 5445.06 & 19124 & 16335.18 \\ \hline
3.0 & 9378 & 10116 & 8179 & 11699 & 10721.16 & 35287 & 32163.49 \\ \hline
\end{tabular}
\caption{Values of $R_c(k,K^+_0,\rho_0,(10^3,10^4))$ for $K^+_0 = \bQ(\zeta_7+\zeta_7^{-1})$.}
\label{real_g3_cyclo7+}

\end{table}